\documentclass[aap,MSNbibl,seceqn,noautosecdot,dvips]{arximspdf}
\usepackage{algorithmic}
\usepackage{stfloats}
\usepackage{graphicx}

\doi{10.1214/13-AAP973} 
\volume{24}
\issue{5}
\pubyear{2014}
\firstpage{2091}
\lastpage{2142}

\makeatletter
\newcommand{\widebar}{\overline}
\fnbelowfloat
\renewcommand{\mid}{|}
\newcommand{\rrvert}{\vert}
\newcommand{\llvert}{\vert}
\newtheorem{teo}{Theorem}
\newproclaim{defn}{Definition}
\newtheorem{prop}{Proposition}
\newtheorem{lem}{Lemma}
\newproclaim{remark}{Remark}
\newproclaim{case}{Case}
\newtheorem{conj}{Conjecture}
\makeatother

\begin{document}
\begin{frontmatter}

\title{Queuing with future information}
\runtitle{Queuing with future information}

\begin{aug}
\author[A]{\fnms{Joel} \snm{Spencer}\ead[label=e1]{spencer@courant.nyu.edu}},
\author[B]{\fnms{Madhu} \snm{Sudan}\ead[label=e2]{madhu@mit.edu}}
\and
\author[C]{\fnms{Kuang} \snm{Xu}\corref{}\ead[label=e3]{kuangxu@mit.edu}\thanksref{t1}}
\runauthor{J. Spencer, M. Sudan and K. Xu}
\affiliation{New York University, Microsoft Research New England and\\
Massachusetts Institute of Technology}
\address[A]{J. Spencer\\
Courant Institute of Mathematical Sciences\\
New York University\\
Room 829, 251 Mercer St.\\
New York, New York 10012\\
USA\\
\printead{e1}} 
\address[B]{M. Sudan\\
Microsoft Research New England\\
One Memorial Drive\\
Cambridge, Massachusetts 02142\\
USA\\
\printead{e2}}
\address[C]{K. Xu\\
Laboratory for Information and Decision Systems (LIDS)\\
Massachusetts Institute of Technology\\
77 Massachusetts Avenue, 32-D666\\
Cambridge, Massachusetts 02139\\
USA\\
\printead{e3}}
\end{aug}
\thankstext{t1}{Supported in part by a research internship at
Microsoft Research New England and by NSF Grants CMMI-0856063 and CMMI-1234062.}

\received{\smonth{1} \syear{2013}}
\revised{\smonth{4} \syear{2013}}

%
\begin{abstract}
We study an admissions control problem, where a queue with service rate
$1-p$ receives incoming jobs at rate $\lambda\in(1-p,1)$, and the
decision maker is allowed to redirect away jobs up to a rate of $p$,
with the objective of minimizing the time-average queue length.

We show that the amount of \emph{information about the future} has a
significant impact on system performance, in the heavy-traffic regime.
When the future is unknown, the optimal average queue length diverges
at rate $\sim\log_{1/(1-p)} \frac{1}{1-\lambda}$, as $\lambda\to
1$. In sharp contrast, when all future arrival and service times are
revealed beforehand, the optimal average queue length converges to a
finite constant, $(1-p)/p$, as $\lambda\to1$. We further show that
the finite limit of $(1-p)/p$ can be achieved using only a \emph
{finite} lookahead window starting from the current time frame, whose
length scales as $\mathcal{O} (\log\frac{1}{1-\lambda}
)$, as $\lambda\to1$. This leads to the conjecture of an interesting
duality between queuing delay and the amount of information about the future.
\end{abstract}

%

\end{frontmatter}

\setcounter{footnote}{1}
\section{\texorpdfstring{Introduction.}{Introduction}}\label{sec1}

\subsection{\texorpdfstring{Variable, but predictable.}{Variable, but predictable}}\label{sec1.1}
The notion of \emph{queues} has been used extensively as a powerful
abstraction in studying dynamic resource allocation systems, where one
aims to match \emph{demands} that arrive over time with available
\emph{resources}, and a queue is used to store currently unprocessed
demands. Two important ingredients often make the design and analysis
of a queueing system difficult: the demands and resources can be both
\emph{variable} and \emph{unpredictable}. \emph{Variability} refers
to the fact that the arrivals of demands or the availability of
resources can be highly volatile and nonuniformly distributed across
the time horizon.\vadjust{\goodbreak} \emph{Unpredictability} means that such
nonuniformity ``tomorrow'' is unknown to the decision maker ``today,''
and she is obliged to make allocation decisions only based on the state
of the system at the moment, and some statistical estimates of the future.

%
\begin{figure}

\includegraphics{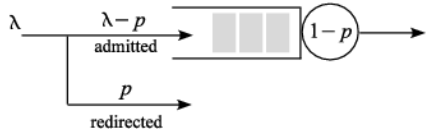}

\caption{An illustration of the admissions control problem, with a
constraint on the a rate of redirection.}\label{figdelill}
\end{figure}

While the world will remain volatile as we know it, in many cases, the
amount of {unpredictability about the future} may be reduced thanks to
\emph{forecasting} technologies and the increasing accessibility of
data. For instance:
\begin{longlist}[(3)]
\item[(1)] advance booking in the hotel and textile industries allows
for accurate forecasting of demands ahead of time \cite{FR96};
\item[(2)] the availability of monitoring data enables traffic
controllers to predict the traffic pattern around potential bottlenecks
\cite{SWO02};
\item[(3)] advance scheduling for elective surgeries could inform care
providers several weeks before the intended appointment \cite{KH02}.
\end{longlist}
In all of these examples, future demands remain \emph{exogenous} and
variable, yet the decision maker is revealed with (some of) their realizations.

\emph{Is there significant performance gain to be harnessed by}
``\emph{looking into the future}?'' In this paper we provide a largely
affirmative answer, in the context of a class of admissions control problems.

\subsection{\texorpdfstring{Admissions control viewed as resource allocation.}{Admissions control viewed as resource allocation}}\label{secintroadminresourcepool} We begin by informally describing our
problem. Consider a single queue equipped with a server that runs at
rate $1-p$ jobs per unit time, where $p$ is a fixed constant in
$(0,1)$, as depicted in Figure~\ref{figdelill}. The queue receives a
stream of incoming jobs, arriving at rate $\lambda\in(0,1)$. If
$\lambda> 1-p$, the arrival rate is greater than the server's
processing rate, and some form of \emph{admissions control} is
necessary in order to keep the system stable. In particular, upon its
arrival to the system, a job will either be \emph{admitted} to the
queue, or \emph{redirected}. In the latter case, the job does not join
the queue, and, from the perspective of the queue, disappears from the
system entirely. The goal of the decision maker is to minimize the
average delay experienced by the admitted jobs, while obeying the
constraint that the average rate at which jobs are redirected \emph
{does not exceeded $p$}.\footnote{Note that as $\lambda\to1$, the
minimum rate of admitted jobs, $\lambda-p$, approaches the server's
capacity $1-p$, and hence we will refer to the system's behavior when
$\lambda\to1$ as the \emph{heavy-traffic regime}.}


One can think of our problem as that of \emph{resource allocation},
where a decision maker tries to match incoming demands with two types
of processing resources: a~\emph{slow local resource} that corresponds
to the server and a \emph{fast external resource} that can process any
job redirected to it almost instantaneously. Both types of resources
are \emph{constrained}, in the sense that their capacities ($1-p$ and
$p$, resp.) cannot change over time, by physical or contractual
predispositions. The processing time of a job at the fast resource is
\emph{negligible compared to that at the slow resource}, as long as
the rate of redirection to the fast resource stays below $p$ in the
long run. Under this interpretation, minimizing the average delay
across \emph{all} jobs is equivalent to minimizing the average delay
across just the \emph{admitted} jobs, since the jobs redirected to the
fast resource can be thought of being processed immediately and
experiencing no delay at all.

For a more concrete example, consider a web service company that enters
a long-term contract with an external cloud computing provider for a
fixed amount of computation resources (e.g., virtual machine instance
time) over the contract period.\footnote{\emph{Example}. As of
September 2012, Microsoft's Windows Azure cloud services offer a
6-month contract for \$71.99 per month, where the client is entitled
for up to 750 hours of virtual machine (VM) instance time each month,
and any additional usage would be charged at a 25\% higher rate. Due to
the large scale of the Azure data warehouses, the speed of any single
VM instance can be treated as roughly constant and independent of the
total number of instances that the client is running concurrently.}
During the contract period, any incoming request can be either served
by the in-house server (slow resource), or be redirected to the cloud
(fast resource), and in the latter case, the job does not experience
congestion delay since the scalability of cloud allows for multiple VM
instance to be running in parallel (and potentially on different
physical machines). The decision maker's constraint is that the total
amount of redirected jobs to the cloud must stay below the amount
prescribed by the contract, which, in our case, translates into a
maximum redirection rate over the contract period. Similar scenarios
can also arise in other domains, where the slow versus fast resources
could, for instance, take on the forms of:
\begin{longlist}[(3)]
\item[(1)] an in-house manufacturing facility versus an external contractor;
\item[(2)] a slow toll booth on the freeway versus a special lane that
lets a car pass without paying the toll;
\item[(3)] hospital bed resources within a single department versus a
cross-departmen\-tal central bed pool.
\end{longlist}

\begin{figure}[t]

\includegraphics{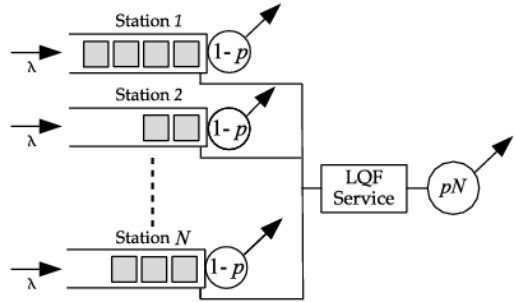}

\caption{Illustration of a model for resource pooling with distributed
and centralized resources~\cite{TX12}.}\label{figpooling}\vspace*{-3pt}
\end{figure}

In a recent work \cite{TX12}, a mathematical model was proposed to
study the benefits of resource pooling in large scale queueing systems,
which is also closely connected to our problem. They consider a
multi-server system where a fraction $1-p$ of a total of $N$ units of
processing resources (e.g., CPUs) is distributed among a set of $N$
local servers,\vadjust{\goodbreak} each running at rate $1-p$, while the remaining fraction
of $p$ is being allocated in a centralized fashion, in the form of a
central server that operates at rate $pN$ (Figure~\ref{figpooling}).
It is not difficult to see, when $N$ is large, the central server
operates at a significantly faster speed than the local servers, so
that a job processed at the central server experiences little or no
delay. In fact, the admissions control problem studied in this paper is
essentially the problem faced by one of the local servers, in the
regime where $N$ is large (Figure~\ref{figpoolingcqueue}). This
connection is explored in greater detail in Appendix~\ref{secresourcepooling}, where we discuss what the implications of our
results in context of resource pooling systems.

%

%
\begin{figure}[b]\vspace*{-3pt}

\includegraphics{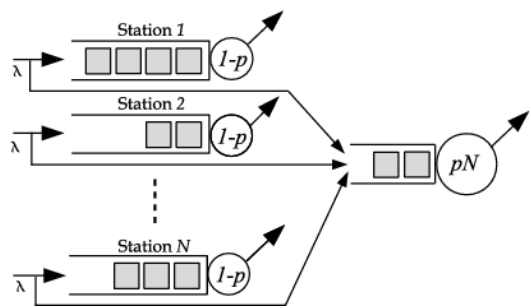}

\caption{Resource pooling using a central queue.}\label{figpoolingcqueue}
\end{figure}

\subsection{\texorpdfstring{Overview of main contributions.}{Overview of main contributions}}\label{sec1.3}
We preview some of the main results in this section. The formal
statements will be given in Section~\ref{secresults}.

\subsubsection{\texorpdfstring{Summary of the problem.}{Summary of the problem}}\label{sec1.3.1}
We consider a continuous-time admissions control problem, depicted in
Figure~\ref{figdelill}. The problem is characterized by three
parameters: $\lambda, p$ and $w$:
\begin{longlist}[(2)]
\item[(1)] Jobs arrive to the system at a rate of $\lambda$ jobs per
unit time, with $\lambda\in(0,1)$. The server operates at a rate of
$1-p$ jobs per unit time, with $p\in(0,1)$.\vadjust{\goodbreak}

\item[(2)] The decision maker is allowed to decide whether an arriving
job is admitted to the queue, or redirected away, with the goal of
minimizing the time-average queue length,\footnote{By Little's law,
the average queue length is essentially the same as average delay, up
to a constant factor; see Section~\ref{secperformancemeasure}.} and
subject to the constraint that the time-average rate of redirection
does not exceed $p$ jobs per unit time.

\item[(3)] The decision maker has access to \emph{information about
the future}, which takes the form of a \emph{lookahead window} of
length $w\in\mathbb{R}_{+}$. In particular, at any time $t$, the
times of arrivals and service availability within the interval
$[t,t+w]$ are revealed to the decision maker. We will consider the
following cases of $w$:
\begin{enumerate}[(a)]
\item[(a)]$w=0$, the \emph{online problem}, where no future information
is available.
\item[(b)]$w=\infty$, the \emph{offline problem}, where entire the future
has been revealed.
\item[(c)]$0<w<\infty$, where future is revealed only up to a finite
lookahead window.
\end{enumerate}
\end{longlist}

Throughout, we will fix $p\in(0,1)$, and be primarily interested in
the system's behavior in the \emph{heavy-traffic regime} of $\lambda
\to1$.

\subsubsection{\texorpdfstring{Overview of main results.}{Overview of main results}}\label{sec1.3.2}

Our main contribution is to demonstrate that the performance of a
redirection policy is highly sensitive to the amount of future
information available, measured by the value of $w$.

\begin{figure}

\includegraphics{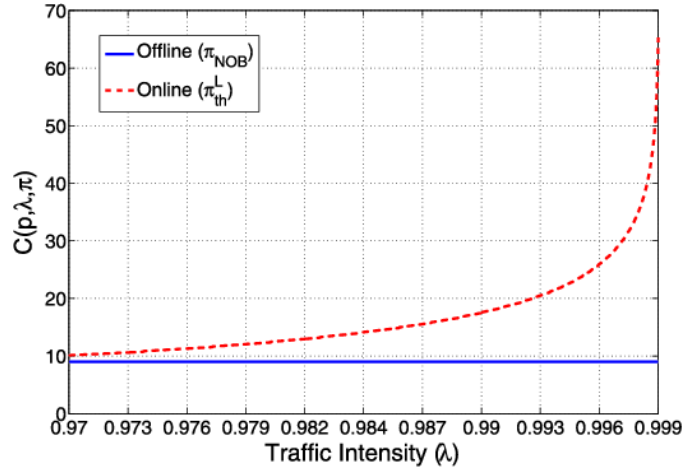}

\caption{Comparison of optimal heavy-traffic delay scaling between
online and offline policies, with $p=0.1$ and $\lambda\to1$. The
value $C(p,\lambda,\pi)$ is the resulting average queue length as a
function of $p$, $\lambda$ and a policy $\pi$.}\label{figdelayscale}
\end{figure}

Fix $p\in(0,1)$, and let the arrival and service processes be Poisson.
For the online problem ($w=0$), we show the optimal time-average queue
length, $C^{\mathrm{opt}}_0$, approaches infinity in the heavy-traffic regime,
at the rate
\[
C^{\mathrm{opt}}_0 \sim\log_{1/(1-p)}\frac{1}{1-\lambda}\qquad\mbox{as }\lambda\to1.
\]
In sharp contrast, the optimal average queue length among offline
policies \mbox{($w=\infty$)}, $C^{\mathrm{opt}}_{\infty}$, converges to a \emph{constant},
\[
C^{\mathrm{opt}}_\infty\to\frac{1-p}{p}\qquad\mbox{as }\lambda\to1
\]
and this limit is achieved by a so-called no-job-left-behind policy.
Figure~\ref{figdelayscale} illustrates this difference in delay
performance for a particular value of $p$.

%

Finally, we show that the no-job-left-behind policy for the offline
problem can be modified, so that the \emph{same} optimal heavy-traffic
limit of $\frac{1-p}{p}$ is achieved even with a \emph{finite}
lookahead window, $w(\lambda)$, where
\[
w(\lambda) = \mathcal{O} \biggl(\log\frac{1}{1-\lambda} \biggr)\qquad\mbox{as }
\lambda\to1.
\]
This is of practical importance because in any realistic application,
only a finite amount of future information can be obtained.

On the methodological end, we use a sample path-based framework to
analyze the performance of the offline and finite lookahead policies,
borrowing tools from renewal theory and the theory of random walks. We
believe that our techniques could be substantially generalized to
incorporate general arrival and service processes, diffusion
approximations as well as observational noises. See Section~\ref{secconclusions} for a more elaborate discussion.

\subsection{\texorpdfstring{Related work.}{Related work}}\label{sec1.4}

There is an extensive body of work devoted to various Markov (or \emph
{online}) admissions control problems; the reader is referred to the
survey of \cite{Sti85} and references therein. Typically, the problem
is formulated as an \mbox{instance} of a Markov decision problem (MDP), where
the decision maker, by \mbox{admitting} or rejecting incoming jobs, seeks to
maximize a long-term average objective \mbox{consisting} of rewards (e.g.,
throughput) minus costs (e.g., waiting time experienced by a customer).
The case where the maximization is performed subject to a constraint on
some average cost has also been studied, and it has been shown, for a
family of reward and cost functions, that an optimal policy assumes a
``threshold-like'' form, where the decision maker redirects the next
job only if the current queue length is great or equal to $L$, with
possible randomization if at level $L-1$, and always admits the job if
below $L-1$; cf.~\cite{BR86}. Indeed, our problem, where one tries to
minimize average queue length (delay) subject to a lower-bound on the
throughput (i.e., a maximum redirection rate), can be shown to belong
to this category, and the online heavy-traffic scaling result is a
straightforward extension following the MDP framework, albeit dealing
with technicalities in extending the threshold characterization to an
infinite state space, since we are interested in the regime of $\lambda
\to1$.

However, the resource allocation interpretation of our admissions
control problem as that of matching jobs with fast and slow resources,
and, in particular, its connections to resource pooling in the
many-server limit, seems to be largely unexplored. The difference in
motivation perhaps explains why the optimal online heavy-traffic delay
scaling of $\log_{1/(1-p)}\frac{1}{1-\lambda}$ that emerges by
fixing $p$ and taking $\lambda\to1$ has not appeared in the
literature, to the best our knowledge.

There is also an extensive literature on \emph{competitive analysis},
which focuses on the \emph{worst-case} performance of an online
algorithms compared to that of an optimal offline version (i.e.,
knowing the entire input sequence). The reader is referred to~\cite
{AE98} for a comprehensive survey, and the references therein on
packing-type problems, such as load balancing and machine scheduling
\cite{Aza98}, and call admission and routing \cite{BAP93}, which are
more related to our problem. While our optimality result for the policy
with a finite lookahead window is stated in terms of the \emph
{average} performance given stochastic inputs, we believe that the
analysis can be extended to yield worst-case competitive ratios under
certain input regularity conditions.

In sharp contrast to our knowledge of the online problems,
significantly less is known for settings in which information about the
future is taken into consideration. In \cite{Naw90}, the author
considers a variant of the flow control problem where the decision
maker knows the job size of the arriving customer, as well as the
arrival and time and job size of the next customer, with the goal of
maximizing certain discounted or average reward. A characterization of
an optimal stationary policy is derived under a standard semi-Markov
decision problem framework, since the lookahead is limited to the next
arriving job. In \cite{CH93}, the authors consider a scheduling
problem with one server and $M$ parallel queues, motivated by
applications in satellite systems where the link qualities between the
server and the queues vary over time. The authors compare the
throughput performance between several online policies with that of an
offline policy, which has access to all future instances of link
qualities. However, the offline policy takes the form of a Viterbi-like
dynamic program, which, while being throughput-optimal by definition,
provides limited qualitative insight.

One challenge that arises as one tries to move beyond the online
setting is that policies with lookahead typically do not admit a clean
Markov description, and hence common techniques for analyzing Markov
decision problems do not easily apply. To circumvent the obstacle, we
will first relax our problem to be fully offline, which turns out to be
surprisingly amenable to analysis. We then use the insights from the
optimal offline policy to construct an optimal policy with a finite
look-ahead window, in a rather straightforward manner.

In other application domains, the idea of exploiting future information
or predictions to improve decision making has been explored. Advance
reservations (a~form of future information) have been studied in lossy
networks \cite{CJP99,LR07} and, more recently, in revenue management
\cite{SL11}. Using simulations, \cite{KH02} demonstrates that the use
of a one-week and two-week advance scheduling window for elective
surgeries can improve the efficiency at the associated intensive care
unit (ICU). The benefits of advanced booking program for supply chains
have been shown in \cite{FR96} in the form of reduced demand
uncertainties. While similar in spirit, the motivations and dynamics in
these models are very different from ours.

Finally, our formulation of the slow an fast resources had been in part
inspired by the literature of resource pooling systems, where one
improves overall system performance by (partially) sharing individual
resources in collective manner. The connection of our problem to a
specific multi-server model proposed by \cite{TX12} is discussed in
Appendix~\ref{secresourcepooling}. For the general topic of resource
pooling, interested readers are referred to \cite{MR98,HL99,BW01,MS04} and the references therein.

\subsection{\texorpdfstring{Organization of the paper.}{Organization of the paper}}\label{sec1.5}

The rest of the paper is organized as follows. The mathematical model
for our problem is described in Section~\ref{secmodel}. Section~\ref{secresults} contains the statements of our main results, and
introduces the no-job-leftb-behind policy ($\pi_{\mathrm{NOB}}$), which will be
a central object of study for this paper. Section~\ref{secinterpret}
presents two alternative descriptions of the no-job-left-behind policy
that have important structural, as well as algorithmic, implications.
Sections~\ref{secoptonline}--\ref{secfinitelookahead} are
devoted to the proofs for the results concerning the online, offline
and finite-lookahead policies, respectively. Finally, Section~\ref{secconclusions} contains some concluding remarks and future directions.

\section{\texorpdfstring{Model and setup.}{Model and setup}}\label{secmodel}

\subsection{\texorpdfstring{Notation.}{Notation}}\label{sec2.1}
We will denote by $\mathbb{N}$, $\mathbb{Z}_{+}$ and $\mathbb
{R}_{+}$, the set of natural numbers, nonnegative integers and
nonnegative reals, respectively. Let $f,g\dvtx  \mathbb{R}_{+}\to\mathbb
{R}_{+}$ be two functions. We will use the following asymptotic
notation throughout: $f(x)\lesssim g(x)$ if $\lim_{x\to1}\frac
{f(x)}{g(x)}\leq1$, $f(x)\gtrsim g(x)$ if $\lim_{x\to1}\frac
{f(x)}{g(x)}\geq1$; $f(x)\ll g(x)$ if $\lim_{x\to1}\frac
{f(x)}{g(x)}= 0$ and $f(x)\gg g(x)$ if $\lim_{x\to1}\frac
{f(x)}{g(x)} = \infty$.

\subsection{\texorpdfstring{System dynamics.}{System dynamics}}\label{sec2.2}
An illustration of the system setup is given in Figure~\ref
{figdelill}. The system consists of a single-server queue running in
continuous time ($t\in\mathbb{R}_+)$, with an unbounded buffer that
stores all unprocessed jobs. The queue is assumed to be empty at $t=0$.

Jobs arrive to the system according to a Poisson
process with rate $\lambda$, $\lambda\in(0,1 )$, so that
the intervals between two adjacent arrivals are independent
and exponentially distributed with mean $\frac{1}{\lambda}$. We will
denote by $ \{A(t)\dvtx  t\in\mathbb{R}_{+} \}$ the cumulative
arrival process,
where $A(t)\in\mathbb{Z}_{+}$ is the total number of arrivals to the
system by time
$t$.

The processing of jobs by the server is modeled
by a Poisson process of rate $1-p$. When the service process
receives a jump at time $t$, we say that a service token is generated.
If the queue is not empty at time $t$, exactly one job ``consumes''
the service token and leaves the system immediately. Otherwise, the
service token is ``wasted'' and has no impact on the future evolution
of the system.\footnote{When the queue is nonempty, the generation of
a token can be interpreted as the completion of a previous job, upon
which the server is ready to fetch the next job. The time between two
consecutive tokens corresponds to the service time. The waste of a
token can be interpreted as the server starting to serve a ``dummy
job.'' Roughly speaking, the service token formulation, compared to
that of a constant speed server processing jobs with exponentially
distributed sizes, provides a performance upper-bound due to the
inefficiency caused by dummy jobs, but has very similar performance in
the heavy-traffic regime, in which the tokens are almost never wasted.
Using such a point process to model services is not new, and the reader
is referred to \cite{TX12} and the references therein.

It is, however, important to note a key assumption implicit in the
service token formulation: the processing times are intrinsic to the
server, and \emph{independent} of the job being processed. For
instance, the sequence of service times will not depend on the order in
which the jobs in the queue are served, so long as the server remains
busy throughout the period. This distinction is of little relevance for
an $M/M/1$ queue, but can be important in our case, where the
redirection decisions may depend on the future. See discussion in
Section~\ref{secconclusions}.} We will denote by $ \{S(t)\dvtx  t\in
\mathbb{R}_{+} \}$ the cumulative token generation
process, where $S(t)\in\mathbb{Z}_{+}$ is the total number of service
tokens generated
by time $t$.

When $\lambda>1-p$, in order to maintain the stability of the queue, a
decision maker has the option of ``redirecting'' a job \emph{at the
moment of its arrival}. Once redirected, a job effectively
``disappears,'' and for this reason, we will use the word \textit{deletion} as a synonymous term for redirection throughout the rest of
the paper, because it is more intuitive to think of deleting a job in
our subsequent sample-path analysis. Finally, the decision maker is
allowed to delete up to a time-average rate of $p$.

\subsection{\texorpdfstring{Initial sample path.}{Initial sample path}}\label{sec2.3}

Let $ \{ Q^{0} (t )\dvtx {t\in\mathbb{R}_{+}} \}$ be
the continuous-time queue length process, where $Q^0(t)\in\mathbb
{Z}_{+}$ is the queue length at time $t$ if \emph{no deletion} is
applied at any time. We say that an \textit{event} occurs at time $t$ if
there is either an arrival, or a generation of service token,
at time $t$. Let $T_{n}$, $n\in\mathbb{N}$, be the time of the $n$th
event in the system.
Denote by $ \{ Q^{0} [n ]\dvtx n\in\mathbb{Z}_{+} \}
$ the embedded
discrete-time process of $ \{Q^{0} ( t ) \}$,
where $Q^0 [n ]$ is the length of the queue sampled immediately
after the $n$th event,\footnote{The notation $f(x-)$ denotes the right-limit of $f$ at $x$:
$f(x-)=\lim_{y\downarrow x}f(y)$.
In this particular context, the values of $Q^{0} [n ]$ are
well defined, since the sample paths of Poisson processes are
right-continuous-with-left-limits
(RCLL) almost surely.}
\[
Q^{0} [n ]= Q^{0} (T_{n}- ),\qquad n\in \mathbb{N}
\]
with the initial condition $Q^0[0]=0$. It is well known that $Q^0$ is a
random walk on~$\mathbb{Z}_{+}$, such that for all $x_1,x_2\in\mathbb
{Z}_{+}$ and $n\in\mathbb{Z}_{+}$,
%
%
\begin{equation}
\qquad \mathbb{P} \bigl(Q^0[n+1]=x_2 \mid Q^0[n]=x_1
\bigr) = \cases{ \displaystyle\frac{\lambda}{\lambda+1-p}, &\quad$x_2-x_1=1$,
\vspace*{5pt}\cr
\displaystyle\frac{1-p}{\lambda+1-p}, &\quad$x_2-x_1=-1$,
\vspace*{5pt}\cr
0, &\quad otherwise,} \label{eqQ0trans1}
\end{equation}
if $x_1>0$ and
%
%
\begin{equation}
\mathbb{P} \bigl(Q^0[n+1]=x_2 \mid Q^0[n]=x_1
\bigr) = \cases{ \displaystyle\frac{\lambda}{\lambda+1-p}, &\quad$x_2-x_1=1$,
\vspace*{5pt}\cr
\displaystyle\frac{1-p}{\lambda+1-p}, &\quad$ x_2-x_1=0$,
\vspace*{5pt}\cr
0, &\quad otherwise,}\label{eqQ0trans2}
\end{equation}
if $x_1=0$. Note that, when $\lambda>1-p$, the random walk $Q^0$ is transient.

The process $Q^0$ contains \emph{all relevant information} in the
arrival and service processes, and will be the main object of study of
this paper. We will refer to $Q^0$ as the \emph{initial sample path}
throughout the paper, to distinguish it from sample paths obtained
after deletions have been made.

\subsection{\texorpdfstring{Deletion policies.}{Deletion policies}}\label{sec2.4}

Since a deletion can only take place when there is an arrival,
it suffices to define the locations of deletions with respect to the
discrete-time process $ \{Q^0[n]\dvtx n\in\mathbb{Z}_{+} \}$,
and throughout,
our analysis will focus on discrete-time queue length processes unless
otherwise specified. Let $\Phi(Q )$ be the locations of
all arrivals
in a discrete-time queue length process $Q$, that is,
\[
\Phi(Q )= \bigl\{ n\in\mathbb{N}\dvtx  Q [n ]>Q [n-1 ] \bigr\}
\]
and for any $M\subset\mathbb{Z}_{+}$, define the counting process
$ \{I(M,n)\dvtx n\in\mathbb{N} \}$ associated with $M$
as\footnote{$|X|$ denotes the cardinality of $X$.}
%
%
\begin{equation}
\label{eqIdef} I(M,n) = \bigl\llvert\{1, \ldots, n \}\cap M \bigr
\rrvert.
\end{equation}

%
\begin{defn}[(Feasible deletion sequence)]\label{deffeasibleDel}
The sequence $M= \{ m_{i} \} $ is said to be a {feasible
deletion sequence} with respect to a discrete-time queue length process,
$Q^{0}$, if all of the following hold:
\begin{longlist}[(2)]
\item[(1)] All elements in $M$ are unique, so that at most one
deletion occurs
at any slot.
\item[(2)] $M\subset\Phi(Q^{0} )$, so that a deletion
occurs only when there
is an arrival.
\item[(3)]
%
%
\begin{equation}
\limsup_{n\rightarrow\infty}\frac{1}{n} I (M,n )\leq\frac{p}{\lambda+(1-p)}\qquad\mbox{a.s.} \label{eqrateconstr}
\end{equation}
so that the time-average deletion rate is at most $p$.
\end{longlist}
In general, $M$ is also allowed to be a finite set.
\end{defn}
The denominator $\lambda+ (1-p )$ in equation~(\ref
{eqrateconstr}) is due to the fact that the total rate of events in
the system is $\lambda+ (1-p )$.\footnote{This is equal to
the total rate of jumps in $A(\cdot)$ and $S(\cdot)$.} Analogously,
the deletion rate in continuous time is defined by
%
%
\begin{equation}
r_d = (\lambda+1-p)\cdot\limsup_{n\rightarrow\infty}
\frac{1}{n} I (M,n ).
\end{equation}

The impact of a deletion sequence to the evolution of the queue length
process is formalized in the following definition.

%
\begin{defn}[(Deletion maps)]\label{defdelMap}
Fix an initial queue length process $ \{
Q^{0}[n]\dvtx\break  n\in\mathbb{N} \}$ and a corresponding feasible
deletion sequence $M= \{ m_{i} \}$.
\begin{longlist}[(2)]
\item[(1)] The \textit{point-wise deletion map} $D_{P}
(Q^{0},m )$
outputs the resulting process after a deletion is made to $Q^{0}$
in slot $m$. Let $Q'=D_{P} (Q^{0},m )$. Then
%
%
\begin{equation}
Q' [n ]=\cases{ Q^{0} [n ]-1, &\quad$n\geq m$ and
$Q^0[t]>0\ \forall t\in\{m,\ldots,n \}$;
\vspace*{5pt}\cr
Q^{0} [n ],
&\quad otherwise,}\label{eqxitr}
\end{equation}

\item[(2)] the \textit{multi-point deletion map }$D (Q^{0},M
)$ outputs
the resulting process after all deletions in the set $M$ are made to
$Q^{0}$. Define $Q^{i}$ recursively as $Q^{i}=D_{P}
(Q^{i-1},m_{i} )$, $\forall i\in\mathbb{N}$. Then $Q^{\infty
}=D (Q^{0},M )$ is defined as the point-wise limit
%
%
\begin{equation}
Q^{\infty}[n]=\lim_{i\rightarrow\min\{|M|,\infty\}
}Q^{i}[n]\qquad\forall n \in\mathbb{Z}_{+}. \label{eqqinft}
\end{equation}
\end{longlist}
\end{defn}

The definition of the point-wise deletion map reflects the earlier
assumption that the service time of a job only depends on the speed of
the server at the moment and is independent of the job's identity; see
Section~\ref{secmodel}. Note also that the value of $Q^{\infty}
[n ]$ depends only
on the total number of deletions before $n$ [equation~(\ref{eqxitr})],
which is at most $n$, and the limit in equation~(\ref{eqqinft}) is
justified. Moreover, it is not difficult to see that the order in which
the deletions are made has no impact on the resulting sample path, as
stated in the lemma below. The proof is omitted.
%
\begin{lem}
\label{lemlocindp} Fix an initial sample path $Q^0$, and let $M$ and
$\widetilde{M}$ be two feasible deletion sequences that contain the same
elements. Then $D (Q^{0},M )=D (Q^{0},\widetilde
{M} )$.
\end{lem}

We next define the notion of a deletion policy that outputs a deletion
sequence based on the (limited) knowledge of an initial sample path
$Q^0$. Informally, a deletion policy is said to be $w$-lookahead if it
makes its deletion decisions based on the knowledge of $Q^0$ up to $w$
units of time into the future (in continuous time).

%
\begin{defn}[($w$-lookahead deletion policies)] \label{defw-pred} Fix $w\in
\mathbb{R}_{+}\cup\{\infty\}$. Let $\mathcal
{F}_t=\sigma(Q^0(s);s\leq t )$ be the natural filtration
induced by $ \{Q^0(t)\dvtx t\in\mathbb{R}_{+} \}$ and
$\mathcal{F}_\infty= \bigcup_{t\in\mathbb{Z}_{+}}\mathcal{F}_t$. A
{$w$-predictive deletion policy} is a mapping, $\pi\dvtx  \mathbb
{Z}_{+}^{\mathbb{R}_+}\rightarrow\mathbb{N}^{\infty}$, such that:
\begin{longlist}[(2)]
\item[(1)] $M=\pi(Q^{0} )$
is a feasible deletion sequence a.s.;
\item[(2)] $ \{n\in M \}$ is $\mathcal{F}_{T_n+w}$
measurable, for all $n\in\mathbb{N}$.
\end{longlist}
We will denote by $\Pi_{w}$ the family of all $w$-lookahead deletion
policies.
\end{defn}

The parameter $w$ in Definition~\ref{defw-pred} captures the amount
of information that the deletion policy has about the future:
\begin{longlist}[(2)]
\item[(1)] When $w=0$, all deletion decisions are made solely based on
the knowledge of the system up to the current time frame. We will refer
to $\Pi_{0}$ as \textit{online policies}.
\item[(2)] When $w=\infty$, the entire sample path of $Q^0$ is
revealed to the decision maker at $t=0$. We will
refer to $\Pi_{\infty}$ as \textit{offline policies}.
\item[(3)] We will refer to $\Pi_w, 0<w<\infty$, as policies with a
\emph{lookahead window of size $w$}.
\end{longlist}

\subsection{\texorpdfstring{Performance measure.}{Performance measure}}\label{secperformancemeasure}
Given a discrete-time queue length process $Q$ and $n\in\mathbb{N}$,
denote by
$S (Q,n )\in\mathbb{Z}_{+}$ the partial sum
%
%
\begin{equation}
S (Q,n )=\sum_{k=1}^{n}Q [k
].\label{eqpartials}
\end{equation}

%
\begin{defn}[(Average post-deletion queue length)] Let $Q^{0}$
be an initial queue length process. Define $C(p,\lambda,\pi) \in
\mathbb{R}_{+}$
as the expected average queue length after applying a deletion policy
$\pi$,
%
%
\begin{equation}
C(p,\lambda,\pi)=\mathbb{E} \biggl(\limsup_{n\rightarrow\infty}
\frac{1}{n}S \bigl(Q_{\pi}^{\infty},n \bigr) \biggr),
\label{eqC}
\end{equation}
where $Q_{\pi}^{\infty}=D (Q^{0},\pi(Q^{0} ) )$,
and the expectation is taken over all realizations of~$Q^{0}$ and the
randomness used by $\pi$ internally, if any.
\end{defn}

\begin{remark*}[(Delay versus queue length)]
By Little's law, the
long-term average waiting time of a typical
customer in the queue is equal to the long-term average queue length
divided by the arrival rate (independent of the service discipline
of the server). Therefore, if our goal is to minimize the average waiting
time of the jobs that remain after deletions, it suffices to use
$C(p,\lambda,\pi)$
as a performance metric in order to judge the effectiveness of a deletion
policy $\pi$. In particular, denote by $T_{\mathrm{all}}\in\mathbb{R}_{+}$
the time-average queueing delay experienced by all jobs, where deleted
jobs are assumed to have a delay of zero, then $\mathbb{E}(T_{\mathrm{all}}) =
\frac{1}{\lambda} C(p,\lambda, \pi)$, and hence the average queue
length and delay coincide in the heavy-traffic regime, as $\lambda\to
1$. With an identical argument, it is easy to see that the average
delay among \emph{admitted} jobs, $T_{\mathrm{adt}}$, satisfies $\mathbb
{E} (T_{\mathrm{adt}} )=\frac{1}{\lambda-r_d} C(p,\lambda, \pi)$,
where $r_d$ is the continuous-time deletion rate under $\pi$.
Therefore, we may use the terms ``delay'' and ``average queue length''
interchangeably in the rest of the paper, with the understanding that
they represent essentially the same quantity up to a constant.

Finally, we define the notion of an optimal delay within a family of policies.
\end{remark*}

%
\begin{defn}[(Optimal delay)] Fix $w\in\mathbb{R}_{+}$. We call $C_{\Pi
_{w}}^{*}(p,\lambda)$
the optimal delay in $\Pi_{w}$, where
%
%
\begin{equation}
C_{\Pi_{w}}^{*}(p,\lambda)=\inf_{\pi\in\Pi_{w}}C(p,
\lambda,\pi).\label{eqLstr}
\end{equation}
\end{defn}

\section{\texorpdfstring{Summary of main results.}{Summary of main results}}\label{secresults}

We state the main results of this paper in this section, whose proofs
will be presented in Sections~\ref{secoptonline}--\ref{secfinitelookahead}.

\subsection{\texorpdfstring{Optimal delay for online policies.}{Optimal delay for online policies}}\label{sec3.1}
%
\begin{defn}[(Threshold policies)]
We say that $\pi_{\mathrm{th}}^{L}$ is an $L$-threshold
policy, if a job arriving at time $t$ is deleted if and only if the
queue length at time $t$ is greater or equal to $L$.
\end{defn}

The following theorem shows that the class of threshold policies
achieves the optimal heavy-traffic delay scaling in $\Pi_0$.

%
\begin{teo}[(Optimal online policies)]\label{teoonline}
Fix $p\in(0,1)$, and let
\[
L (p,\lambda) = \biggl\lceil\log_{\lambda/(1-p)}\frac{p}{1-\lambda} \biggr\rceil.
\]
Then:
\begin{longlist}[(2)]
\item[(1)] $\pi_{\mathrm{th}}^{L (p,\lambda)}$ is feasible for
all $\lambda\in(1-p,1 )$.\vspace*{1pt}
\item[(2)] $\pi_{\mathrm{th}}^{L (p,\lambda)}$ is asymptotically optimal
in $\Pi_{0}$ as $\lambda\to1$,
\[
C \bigl(p,\lambda,\pi_{\mathrm{th}}^{L (p,\lambda)} \bigr)\sim
C_{\Pi_{0}}^{*} (p,\lambda)\sim\log_{1/(1-p)}
\frac
{1}{1-\lambda}\qquad\mbox{as }\lambda\rightarrow1.
\]
\end{longlist}
\end{teo}
\begin{pf}
See Section~\ref{secoptonline}.
\end{pf}

\subsection{\texorpdfstring{Optimal delay for offline policies.}{Optimal delay for offline policies}}\label{sec3.2}

Given the sample path of a random walk $Q$, let $U (Q,n )$
the number of slots till $Q$ reaches the level $Q[n]-1$ after slot $n$:
%
%
\begin{equation}
U (Q,n )=\inf\bigl\{j\geq1\dvtx Q [n+j ]=Q[n]-1 \bigr\}.
\end{equation}

%
\begin{defn}[(No-job-left-behind policy\footnote{The reason for choosing
this name will be made in clear in Section~\ref{secinterpstack},
using the ``stack'' interpretation of this policy.})]\label{defnob}
Given an initial
sample path $Q^{0}$, the no-job-left-behind policy, denoted by $\pi
_{\mathrm{NOB}}$, deletes all arrivals in the set $\Psi$, where
%
%
\begin{equation}
\Psi= \bigl\{n\in\Phi\bigl(Q^0 \bigr)\dvtx  U \bigl(Q^0,n
\bigr)=\infty\bigr\}. \label{eqpsi}
\end{equation}
We will refer to the deletion sequence generated by $\pi_{\mathrm{NOB}}$ as
$M^\Psi= \{m^\Psi_i\dvtx\break  i\in\mathbb{N} \}$, where
$M^\Psi= \Psi$.
\end{defn}

\begin{figure}

\includegraphics{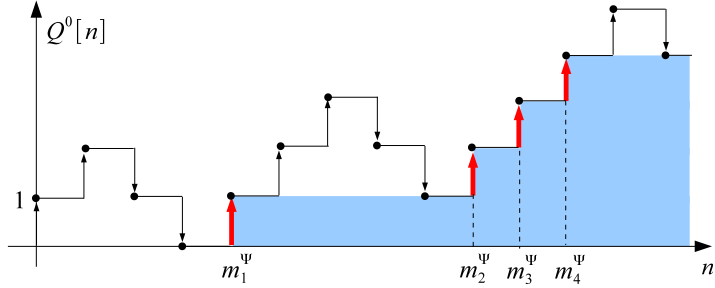}

\caption{Illustration of applying $\pi_{\mathrm{NOB}}$ to an initial sample
path, $Q^0$, where the deletions are marked by bold red arrows.}\label{figwater1}
\end{figure}

\begin{figure}

\includegraphics{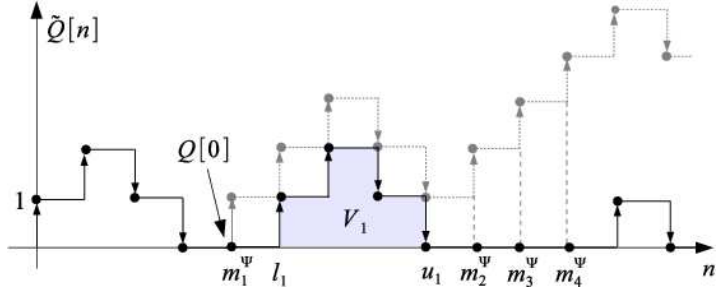}

\caption{The solid lines depict the resulting sample path, $\widetilde{Q}
=D (Q^0,M^\Psi)$, after applying $\pi_{\mathrm{NOB}}$ to
$Q^0$.}\label{figwater2}
\end{figure}

\begin{figure}

\includegraphics{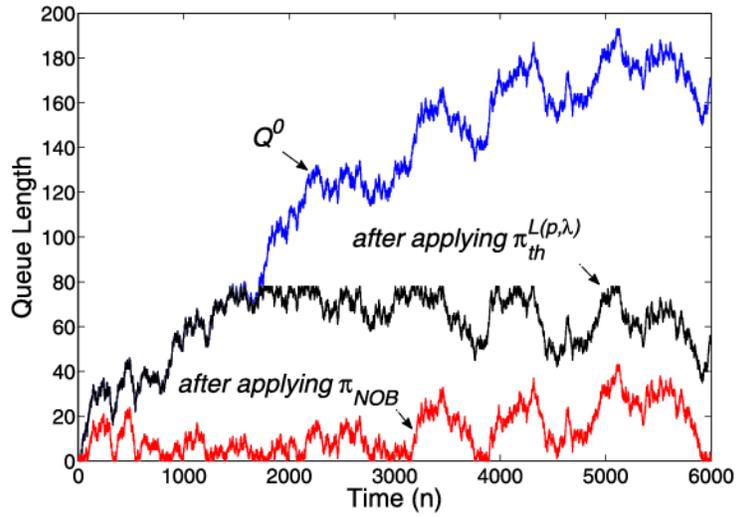}

\caption{Example sample paths of $Q^0$ and those obtained after
applying $\pi_{\mathrm{th}}^{L (p,\lambda)}$ and $\pi_{\mathrm{NOB}}$ to
$Q^0$, with $p=0.05$ and $\lambda=0.999$.}\label{figsamplepaths}
\end{figure}

In other words, $\pi_{\mathrm{NOB}}$ would delete a job arriving at time $t$ if
and only if the initial queue length process never returns to below the
current level in the future, which also implies that
%
%
\begin{equation}
Q^0[n] \geq Q^0 \bigl[m^\Psi_i
\bigr]\qquad\forall i\in\mathbb{N}, n\geq m^\Psi_i.
\label{eqmiinc}
\end{equation}
Examples of the $\pi_{\mathrm{NOB}}$ policy being applied to a particular
sample path are given in Figures~\ref{figwater1} and~\ref{figwater2} (illustration), as well as in Figure~\ref{figsamplepaths} (simulation).

%

%

It turns out that the delay performance of $\pi_{\mathrm{NOB}}$ is about as
good as we can hope for in heavy traffic, as is formalized in the next theorem.

%
\begin{teo}[(Optimal offline policies)]
\label{teooffline}
Fix $p\in(0,1 )$.
\begin{longlist}[(2)]
\item[(1)] $\pi_{\mathrm{NOB}}$ is feasible for all $\lambda\in
(1-p,1 )$ and\footnote{It is easy to see that $\pi_{\mathrm{NOB}}$ is
not a very efficient deletion policy for relatively small values of
$\lambda$. In fact, $C (p,\lambda,\pi_{\mathrm{NOB}} )$ is a \emph
{decreasing} function of $\lambda$. This problem can be fixed by
injecting into the arrival process an Poisson process of ``dummy jobs''
of rate $1-\lambda-\varepsilon$, so that the total rate of arrival is
$1-\varepsilon$, where $\varepsilon\approx0$. This reasoning implies that
$(1-p)/p$ is a uniform upper-bound of $C^*_{\Pi_\infty}(p,\lambda)$
for all $\lambda\in(0,1)$. }
%
%
\begin{equation}
C (p,\lambda,\pi_{\mathrm{NOB}} )=\frac{1-p}{\lambda-(1-p)}.
\end{equation}
\item[(2)] $\pi_{\mathrm{NOB}}$ is asymptotically optimal in $\Pi_{\infty}$
as $\lambda\to1$,
\[
\lim_{\lambda\rightarrow1}C (p,\lambda,\pi_{\mathrm{NOB}} )=\lim
_{\lambda\rightarrow1}C_{\Pi_{\infty}}^{*} (p,\lambda)=
\frac{1-p}{p}.
\]
\end{longlist}
\end{teo}

\begin{pf}
See Section~\ref{secoffline}.
\end{pf}

%
\begin{remark}[(Heavy-traffic ``delay collapse'')]
It
is perhaps surprising to observe that the heavy-traffic scaling
essentially ``collapses'' under $\pi_{\mathrm{NOB}}$: the average queue length
converges to a finite value, $\frac{1-p}{p}$, as $\lambda\to1$,
which is in sharp contrast~with the optimal scaling of $\sim\log
_{1/(1-p)}\frac{1}{1-\lambda}$ for the\vspace*{1pt} online policies, given by
Theorem~\ref{teoonline}; see Figure~\ref{figdelayscale} for an
illustration of this difference. A ``stack'' interpretation of the
no-job-left-behind policy (Section~\ref{secinterpantireac}) will
help us understand intuitively {why} such a drastic discrepancy exists
between the online and offline heavy-traffic scaling behaviors.

Also, as a by-product of Theorem~\ref{teooffline}, observe that the
heavy-traffic limit scales, in $p$, as
%
%
\begin{equation}
\lim_{\lambda\rightarrow1}C_{\Pi_{\infty}}^{*} (p,\lambda) \sim
\frac{1}{p}\qquad\mbox{as $p \to0$.}
\end{equation}
This is consistent with an intuitive notion of ``flexibility'': delay
should degenerate as the system's ability to redirect away jobs diminishes.
\end{remark}

%
\begin{remark}[(Connections to branching processes and Erd\H {o}s--R\'enyi random graphs)]
Let $d<1<c$ satisfy $de^{-d}=ce^{-c}$. Consider a Galton--Watson
birth process in which each node has $Z$ children, where $Z$ is Poisson
with mean $c$. Conditioning on the finiteness of the process gives
a Galton--Watson process where $Z$ is Poisson with mean $d$. This
occurs in the classical analysis of the Erd\H{o}s--R\'enyi random graph
$G(n,p)$ with $p=c/n$. There will be a giant component and the
deletion of that component gives a random graph $G(m,q)$ with $q=d/m$.
As a rough analogy, $\pi_{\mathrm{NOB}}$ deletes those nodes that would be
in the giant component.
\end{remark}


\subsection{\texorpdfstring{Policies with a finite lookahead window.}{Policies with a finite lookahead window}}\label{sec3.3}

In practice, infinite prediction into the future is certainly too much
to ask for. In this section, we show that a natural modification of
$\pi_{\mathrm{NOB}}$ allows for the \emph{same delay} to be achieved, using
only a \emph{finite} lookahead window, whose length, $w(\lambda)$,
increases to infinity as $\lambda\to1$.\footnote{In a way, this is
not entirely surprising, since the $\pi_{\mathrm{NOB}}$ leads to a deletion
rate of $\lambda-(1-p)$, and there is an additional $p-[\lambda
-(1-p)]=1-\lambda$ unused deletion rate that can be exploited. }

Denote by $w\in\mathbb{R}_{+}$ the size of the lookahead window in
continuous time, and $W(n)\in\mathbb{Z}_{+}$ the window size in the
discrete-time embedded process $Q^0$, starting from slot $n$. Letting
$T_n$ be the time of the $n$th event in the system, then
%
%
\begin{equation}
W(n) = \sup\{ k \in\mathbb{Z}_{+}\dvtx  T_{n+k} \leq
T_{n}+w \}. \label{eqwdiscrete}
\end{equation}

For $x \in\mathbb{N}$, define the set of indices
%
%
\begin{equation}
U (Q,n,x )=\inf\bigl\{j\in\{1,\ldots,x \}\dvtx Q [n+j ]=Q[n]-1 \bigr\}.
\end{equation}

%
\begin{defn}[($w$-no-job-left-behind policy)]
Given an initial sample path
$Q^{0}$ and $w>0$, the $w$-no-job-left-behind policy, denoted by $\pi
_{\mathrm{NOB}}^w$, deletes all arrivals in the set $\Psi^{w}$, where
\[
\Psi^{w} = \bigl\{n\in\Phi\bigl(Q^0 \bigr)\dvtx  U
\bigl(Q^0,n,W(n) \bigr)=\infty\bigr\}.
\]
\end{defn}

It is easy to see that $\pi_{\mathrm{NOB}}^w$ is simply $\pi_{\mathrm{NOB}}$ applied
within the confinement of a finite window: a job at $t$ is deleted if
and only if the initial queue length process does not return to below
the current level \emph{within the next $w$ units of time}, assuming
no further deletions are made. Since the window is finite, it is clear
that $\Psi^w \supset\Psi$ for any $w<\infty$, and hence $C
(p,\lambda,\pi_{\mathrm{NOB}}^w )\leq C (p,\lambda,\pi
_{\mathrm{NOB}} )$ for all $\lambda\in(1-p)$. The only issue now becomes
that of feasibility: by making decision only based on a finite
lookahead window, we may end up deleting at a rate greater than $p$.

The following theorem summarizes the above observations and gives an
upper bound on the appropriate window size, $w$, as a function of
$\lambda$.\footnote{Note that Theorem~\ref{teolookahead} implies
Theorem~\ref{teooffline} and is hence stronger.}

%
\begin{teo}[(Optimal delay scaling with finite lookahead)]\label{teolookahead}
Fix $p\in
(0,1 )$. There exists $C>0$, such that if
\[
w (\lambda)= C\cdot\log\frac{1}{1-\lambda},
\]
then $\pi_{\mathrm{NOB}}^{w(\lambda)}$ is feasible and
%
%
\begin{equation}
C \bigl(p,\lambda,\pi_{\mathrm{NOB}}^{w(\lambda)} \bigr)\leq C (p,\lambda,
\pi_{\mathrm{NOB}} ) = \frac{1-p}{\lambda-(1-p)}.
\end{equation}
Since $C_{\Pi_{w (\lambda)}}^* (p,\lambda)
\geq C_{\Pi_{\infty}}^* (p,\lambda)$ and $C^*_{\Pi
_{w (\lambda)}}(p,\lambda)\leq C (p,\lambda,\pi
_{\mathrm{NOB}}^{w(\lambda)} )$, we also have that
%
%
\begin{equation}
\lim_{\lambda\rightarrow1}C_{\Pi_{w (\lambda
)}}^{*} (p,\lambda) = \lim
_{\lambda\rightarrow1}C_{\Pi
_{\infty}}^{*} (p,\lambda)=
\frac{1-p}{p}.
\end{equation}
\end{teo}

\begin{pf}
See Section~\ref{secpfthmlookahead}.
\end{pf}

\subsubsection{\texorpdfstring{Delay-information duality.}{Delay-information duality}}\label{sec3.3.1}

Theorem~\ref{teolookahead} says that one can attain the same
heavy-traffic delay performance as the optimal offline algorithm if
the size of the lookahead window scales as $\mathcal{O}(\log\frac
{1}{1-\lambda})$. Is this the minimum amount of future information
necessary to achieve the same (or comparable) heavy-traffic delay limit
as the optimal offline policy? We conjecture that this is the case, in
the sense that there exists a matching lower bound, as
follows.\looseness=-1

%
\begin{figure}[b]

\includegraphics{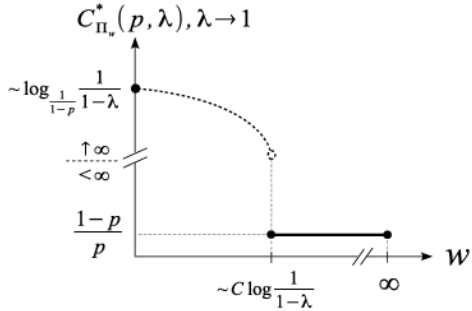}

\caption{``Delay vs. Information.'' Best achievable heavy traffic
delay scaling as a function of the size of the lookahead window, $w$.
Results presented in this paper are illustrated in the solid lines and
circles, and the gray dotted line depicts our conjecture of the unknown
regime of $0<w(\lambda)\lesssim\log(\frac{1}{1-\lambda
} )$.}\label{figcollapse}
\end{figure}

%
\begin{conj}\label{conjinfolowerbound}
Fix $p\in(0,1)$. If $w(\lambda)\ll\log\frac{1}{1-\lambda}$ as
$\lambda\to1$, then
\[
\limsup_{\lambda\to1}C^*_{\Pi_{w(\lambda)}}(p,\lambda) = \infty.
\]
In other words, ``delay collapse'' can occur only if $w(\lambda) =
\Theta(\log\frac{1}{1-\lambda} )$.
\end{conj}

If the conjecture is proven, it would imply a \emph{sharp transition}
in the system's heavy-traffic delay scaling behavior, around the
critical ``threshold'' of $w(\lambda) = \Theta(\log\frac
{1}{1-\lambda} )$. It would also imply the existence of a
symmetric {dual relationship} between \emph{future information} and
\emph{queueing delay}: $\Theta(\log\frac{1}{1-\lambda}
)$ amount of information is required to achieve a finite delay limit,
and one has to suffer $\Theta(\log\frac{1}{1-\lambda}
)$ in delay, if only finite amount of future information is available.

Figure~\ref{figcollapse} summarizes the main results of this paper
from the angle of the delay-information duality. The dotted
line\vadjust{\goodbreak}
segment marks the unknown regime and the sharp transition at its right
endpoint reflects the view of Conjecture~\ref{conjinfolowerbound}.

\section{\texorpdfstring{Interpretations of $\pi_{\mathrm{NOB}}$.}{Interpretations of pi NOB}}\label{secinterpret}
We present two equivalent ways of describing the no-job-left-behind
policy $\pi_{\mathrm{NOB}}$. The \emph{stack interpretation} helps us derive
asymptotic deletion rate of $\pi_{\mathrm{NOB}}$ in a simple manner, and
illustrates the superiority of $\pi_{\mathrm{NOB}}$ compared to an online
policy. Another description of $\pi_{\mathrm{NOB}}$ using time-reversal shows
us that the set of deletions made by $\pi_{\mathrm{NOB}}$ can be calculated
efficiently in linear time (with respect to the length of the time horizon).

\subsection{\texorpdfstring{Stack interpretation.}{Stack interpretation}}\label{secinterpstack}
Suppose that the service discipline adopted by the server is that of
last-in-first-out (LIFO), where the it always fetches a task that
has arrived the latest. In other words, the queue works as a \emph
{stack}. Suppose that we first simulate the stack
without any deletion. It is easy to see that,\vadjust{\goodbreak} when the arrival rate
$\lambda$ is greater than the service rate $1-p$,
there will be a growing set of jobs at the bottom of the stack that
will \emph{never} be processed. Label all such jobs as
``left-behind.'' For example, Figure~\ref{figwater1} shows the
evolution of the queue over time, where all ``left-behind'' jobs are
colored with a blue shade. One can then verify that the policy $\pi
_{\mathrm{NOB}}$ given in Definition~\ref{defnob} is equivalent to deleting
all jobs that are labeled ``left-behind,'' hence the namesake ``No Job
Left Behind.'' Figure~\ref{figwater2} illustrates applying $\pi
_{\mathrm{NOB}}$ to a sample path of $Q^0$, where the $i$th job to be deleted is
precisely the $i$th job among all jobs that would have never been
processed by the server under a LIFO policy.

One advantage of the stack interpretation is that it makes obvious the
fact that the deletion rate induced by $\pi_{\mathrm{NOB}}$ is equal to
$\lambda-(1-p)<p$, as illustrated in the following lemma.

%
\begin{lem} For all $\lambda>1-p$, the following statements hold:
\label{lemnobbasic}
\begin{longlist}[(2)]
\item[(1)] With probability one, there exists $T<\infty$, such that
every service token generated after time $T$ is matched with some job.
In other words, the server never idles after some finite time.
\item[(2)] Let $Q=D (Q^0,M^\Psi)$. We have
%
%
\begin{equation}
\limsup_{n\to\infty}\frac{1}{n}I \bigl(M^\Psi,n
\bigr) \leq\frac{\lambda- (1-p )}{\lambda+1-p}\qquad\mbox{a.s.}, \label{eqIMnlim}
\end{equation}
which implies that $\pi_{\mathrm{NOB}}$ is feasible for all $p\in(0,1)$ and
$\lambda\in(1-p,1)$.
\end{longlist}
\end{lem}

\begin{pf}
See Appendix~\ref{applemnobbasic}.
\end{pf}

\subsubsection{\texorpdfstring{``Anticipation'' vs. ``reaction.''}{``Anticipation'' vs. ``reaction''}}\label{secinterpantireac} Some geometric intuition from the stack
interpretation shows that the power of $\pi_{\mathrm{NOB}}$ essentially
stems
from being highly \emph{anticipatory}. Looking at\vadjust{\goodbreak} Figure~\ref
{figwater1}, one sees that the jobs that are ``left behind'' at the
bottom of the stack correspond to those who arrive during the intervals
where the initial sample path $Q^0$ is taking a consecutive ``upward
hike.'' In other words, $\pi_{\mathrm{NOB}}$ begins to delete jobs when it
anticipates that the arrivals are \emph{just about to} get intense.
Similarly, a job in the stack will be ``served'' if $Q^0$ curves down
eventually in the future, which corresponds $\pi_{\mathrm{NOB}}$'s stopping
deleting jobs as soon as it anticipates that the next few arrivals can
be handled by the server alone. In sharp contrast is the nature of the
optimal online policy, $\pi_{\mathrm{th}}^{L (p,\lambda)}$, which
is by definition ``reactionary'' and begins to delete only when the
current queue length has already reached a high level. The differences
in the resulting sample paths are illustrated via simulations in
Figure~\ref{figsamplepaths}. For example, as $Q^0$ continues to
increase during the first 1000 time slots, $\pi_{\mathrm{NOB}}$ begins deleting
immediately after $t=0$, while no deletion is made by $\pi
_{\mathrm{th}}^{L (p,\lambda)}$ during this period.\looseness=1

As a rough analogy, the offline policy starts to delete \emph{before}
the arrivals get busy, but the online policy can only delete \emph
{after} the burst in arrival traffic has been realized, by which point
it is already ``too late'' to fully contain the delay. This explains,
to certain extend, why $\pi_{\mathrm{NOB}}$ is capable of achieving ``delay
collapse'' in the heavy-traffic regime (i.e., a finite limit of delay
as $\lambda\to1$, Theorem~\ref{teooffline}), while the delay under
even the best online policy diverges to infinity as $\lambda\to1$
(Theorem~\ref{teoonline}).

%
%

\subsection{\texorpdfstring{A linear-time algorithm for $\pi_{\mathrm{NOB}}$.}{A linear-time algorithm for pi NOB}}\label{seclinearalgo}
While the offline deletion problem serves as a nice abstraction, it is
impossible to actually store information about the \emph{infinite}
future in practice, even if such information is available. A natural
finite-horizon version of the offline deletion problem can be posed as
follows: given the values of $Q^0$ over the first $N$ slots, where $N$
finite, one would like to compute the set of deletions made by $\pi_{\mathrm{NOB}}$,
%
%
\[
M^\Psi_N = M^\Psi\cap\{1,\ldots, N\}
\]
assuming that $Q^0[n]> Q^0[N]$ for all $n\geq N$. Note that this
problem also arises in computing the sites of deletions for the $\pi
_{\mathrm{NOB}}^w$ policy, where one would replace $N$ with the length of the
lookahead window, $w$.

We have the following algorithm, which identifies all slots on which a
new ``minimum'' (denoted by the variable $S$) is achieved in $Q^0$,
when viewed in the \emph{reverse} order of time.


It is easy to see that the running time of the above algorithm scales
linearly with respect to the length of the time horizon, $N$. Note that
this is not the unique linear-time algorithm. In fact, one can verify
that the simulation procedure used in describing the stack
interpretation of $\pi_{\mathrm{NOB}}$ (Section~\ref{secinterpret}), which
keeps track of which jobs would eventually be served, is itself a
linear-time algorithm. However, the time-reverse version given here is
arguably more intuitive and simpler to describe.
\eject

\hrule\vspace{0.05in}
{\bf\emph{A linear-time algorithm for $\pi_{\mathrm{NOB}}$}}
\hrule\vspace{0.05in}
\begin{algorithmic}
\STATE$S\leftarrow Q^0[N]$ and $M^\Psi_N\leftarrow\varnothing$
\FOR{$n = N$ down to $1$}
\IF{$Q^0[n]<S$}
\STATE$M^\Psi_N\leftarrow M^\Psi_N\cup\{n+1\}$
\STATE$S \leftarrow Q^0[n]$
\ELSE
\STATE$M^\Psi_N\leftarrow M^\Psi_N$
\ENDIF
\ENDFOR
\RETURN$M^\Psi_N$
\hrule\vspace{0.05in}
\end{algorithmic}

\section{\texorpdfstring{Optimal online policies.}{Optimal online policies}}\label{secoptonline}
Starting from this section and through Section~\ref{secfinitelookahead}, we present the proofs of the results stated in
Section~\ref{secresults}.

We begin with showing Theorem~\ref{teoonline}, by formulating the
online problem as a Markov decision problem (MDP) with an average cost
constraint, which then enables us to use existing results to
characterize the form of optimal policies. Once the family of threshold
policies has been shown to achieve the optimal delay scaling in $\Pi
_0$ under heavy traffic, the exact form of the scaling can be obtained
in a fairly straightforward manner from the steady-state distribution
of a truncated birth--death process.

\subsection{\texorpdfstring{A Markov decision problem formulation.}{A Markov decision problem formulation}}\label{sec5.1}
Since both the arrival and service processes are Poisson, we can
formulate the problem of finding an optimal policy in $\Pi_0$ as a
continuous-time Markov decision problem with an average-cost
constraint, as follows. Let $ \{Q(t)\dvtx t\in\mathbb{R}_{+} \}
$ be the resulting continuous-time queue length process after applying
some policy in $\Pi_0$ to $Q^0$. Let $T_k$ be the $k$th upward jump in
$Q$ and $\tau_k$ the length of the $k$th inter-jump interval, $\tau_k
= T_{k}-T_{k-1}$. The task of a deletion policy, $\pi\in\Pi_0$,
amounts to choosing, for each of the inter-jump intervals, a \emph
{deletion action}, $a_k \in[0,1]$, where the value of $a_k$
corresponds to the probability that the next arrival during the current
inter-jump interval will be deleted. Define $R$ and $K$ to be the \emph
{reward} and \emph{cost} functions of an inter-jump interval, respectively,
%
%
\begin{eqnarray}
\label{eqreward} R(Q_k,a_k,\tau_k) &=&
-Q_k\cdot\tau_k,
\\
K(Q_k,a_k,\tau_k) &=& \lambda(1-a_k)
\tau_k, \label{eqcost}
\end{eqnarray}
where $Q_k = Q(T_k)$. The corresponding MDP seeks to maximize the
time-average reward
%
%
\begin{equation}
\widebar{R}_\pi= \liminf_{n \to\infty} \frac{\mathbb{E}_\pi
(\sum_{k=1}^n R(Q_k,a_k,\tau_k) )}{\mathbb{E}_\pi(\sum_{k=1}^n \tau_k )}
\end{equation}
while obeying the average-cost constraint
%
%
\begin{equation}
\widebar{C}_\pi= \limsup_{n \to\infty} \frac{\mathbb{E}_\pi
(\sum_{k=1}^n K(Q_k,a_k,\tau_k) )}{\mathbb{E}_\pi(\sum_{k=1}^n \tau_k )}
\leq p. \label{eqavgcost}
\end{equation}
To see why this MDP solves our deletion problem, observe that $\widebar
{R}_\pi$ is the negative of the time-average queue length, and
$\widebar
{C}_\pi$ is the time-average deletion rate.

It is well known that the type of constrained MDP described above
admits an optimal policy that is stationary \cite{AS91}, which means
that the action $a_k$ depends solely on the current state, $Q_k$, and
is independent of the time index $k$. Therefore, it suffices to
describe $\pi$ using a sequence, $ \{b_q\dvtx q\in\mathbb
{Z}_{+} \}$, such that $ a_k = b_q$ whenever $Q_k=q$. Moreover,
when the state space is finite,\footnote{This corresponds to a finite
buffer size in our problem, where one can assume that the next arrival
is automatically deleted when the buffer is full, independent of the
value of $a_k$.} stronger characterizations of the $b_q$'s have been
obtained for a family of reward and cost functions under certain
regularity assumptions (Hypotheses 2.7, 3.1 and 4.1 in \cite{BR86}),
which ours do satisfy [equations~(\ref{eqreward}) and (\ref
{eqcost})]. Theorem~\ref{teoonline} will be proved using the next-known result (adapted from Theorem 4.4 in \cite{BR86}):

%
\begin{lem}
\label{lemquasithresh}
Fix $p$ and $\lambda$, and let the buffer size $B$ be finite. There
exists an optimal stationary policy, $ \{b^*_q \}$, of the form
\[
b^*_q=\cases{ 1, &\quad$q< L^*-1$,
\vspace*{3pt}\cr
\xi, &\quad$q=L^*-1$,
\vspace*{3pt}\cr
0, &
\quad$q\geq L^*$}
\]
for some $L^*\in\mathbb{Z}_{+}$ and $\xi\in[0,1]$.
\end{lem}

\subsection{\texorpdfstring{Proof of Theorem \protect\ref{teoonline}.}{Proof of Theorem 1}}\label{secpfthmonline}
\begin{figure}

\includegraphics{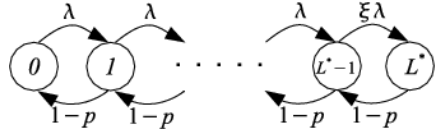}

\caption{The truncated birth--death process induced by $\pi
^*_p$.}\label{figbdchain}
\end{figure}

In words, Lemma~\ref{lemquasithresh} states that the optimal policy admits a
``quasi-threshold'' form: it deletes the next arrival when $Q(t)\geq
L^*$, admits when $Q(t)<L^*-1$, and admits with probability $\xi$ when
$Q(t)=L^*-1$. Suppose, for the moment, that the statements of Lemma
\ref{lemquasithresh} also hold when the buffer size is infinite, an
assumption to be justified by the end of the proof. Denoting by $\pi
^*_p$ the stationary optimal policy associated with $ \{
b^*_q \}$, when the constraint on the average of deletion is $p$
[equation~(\ref{eqavgcost})]. The evolution of $Q(t)$ under $\pi
^*_p$ is that of a birth--death process truncated at state $L^*$, with
the transition rates given in Figure~\ref{figbdchain}, and the
time-average queue length is equal to the expected queue length in
steady state. Using standard calculations involving the steady-state
distribution of the induced Markov process, it is not difficult to
verify that
%
%
\begin{equation}
C\bigl(p,\lambda,\pi_{\mathrm{th}}^{L^* -1}\bigr)\leq C\bigl(p,\lambda,
\pi^*_p\bigr) \leq C\bigl(p,\lambda,\pi_{\mathrm{th}}^{L^*}
\bigr), \label{eqCCC1}
\end{equation}
where $L^*$ is defined as in Lemma~\ref{lemquasithresh}, and
$C(p,\lambda,\pi)$ is the time-average queue length under policy $\pi
$, defined in equation~(\ref{eqC}).

%
%

Denote by $ \{\mu^L_i\dvtx i\in\mathbb{N} \}$ the steady-state
probability of the queue length being equal to $i$, under a threshold
policy $\pi_{\mathrm{th}}^{L}$. Assuming $\lambda\neq1-p$, standard
calculations using the balancing equations yield
%
%
\begin{equation}
\mu^L_i = \biggl(\frac{\lambda}{1-p}\biggr)^i\cdot\biggl(\frac
{1-(\lambda/(1-p))}{1- (\lambda/(1-p))^{L+1}} \biggr)\qquad\forall1\leq i
\leq L \label{eqmui}
\end{equation}
and $\mu_i^L=0$ for all $i \geq L+1$. The time-average queue length is
given by
%
%
\begin{eqnarray}
C\bigl(p,\lambda, \pi_{\mathrm{th}}^{L}\bigr) &=& \sum
_{i=1}^L i\cdot\mu^L_i
\nonumber
\\[-8pt]
\\[-8pt]
&=&\frac{\theta}{ (\theta-1 ) (\theta^{L+1}-1
)}\cdot\bigl[1-\theta^L+L\theta^L(
\theta-1) \bigr],
\nonumber
\end{eqnarray}
where $\theta= \frac{\lambda}{1-p}$. Note that when $\lambda>1-p$,
$\mu_i^L$ is decreasing with respect to $L$ for all $i\in\{
0,1,\ldots,L \}$ [equation~(\ref{eqmui})], which implies that
the time-average queue length is monotonically increasing in $L$, that is,
%
%
\begin{eqnarray}
\label{eqCmonL} &&C\bigl(p,\lambda, \pi_{\mathrm{th}}^{L+1}\bigr) - C
\bigl(p,\lambda, \pi_{\mathrm{th}}^{L}\bigr)
\nonumber
\\
&&\qquad= (L+1)\cdot\mu^{L+1}_{L+1} + \sum
_{i=0}^{L} i\cdot\bigl(\mu^{L+1}_i
- \mu^L_i\bigr)
\nonumber
\\
&&\qquad\geq(L+1)\cdot\mu^{L+1}_{L+1} + L \cdot\Biggl(\sum
_{i=0}^{L} \mu^{L+1}_i
- \mu^L_i \Biggr)
\\
&&\qquad= (L+1)\cdot\mu^{L+1}_{L+1} + L\cdot\bigl(1-\mu
_{i}^{L+1}-1 \bigr)
\nonumber
\\
&&\qquad= \mu^{L+1}_{L+1}>0.
\nonumber
\end{eqnarray}
It is also easy to see that, fixing $p$, since we have that $\theta
>1+\delta$ for all $\lambda$ sufficiently close to $1$, where $\delta
>0$ is a fixed constant, we have
%
%
\begin{equation}
\qquad C\bigl(p,\lambda, \pi_{\mathrm{th}}^{L}\bigr) = \biggl(
\frac{\theta^{L+1}}{\theta
^{L+1}-1} \biggr)L-\frac{\theta}{\theta-1}\cdot\frac{\theta
^{L}-1}{\theta^{L+1}-1} \sim L\qquad\mbox{as $L\to\infty$.} \label{eqCthresh}
\end{equation}

Since deletions only occur when $Q(t)$ is in state $L$, from
equation~(\ref{eqmui}), the average rate of deletions in continuous
time under $\pi_{\mathrm{th}}^{L}$ is given by
%
%
\begin{equation}
\quad r_d \bigl(p,\lambda,\pi_{\mathrm{th}}^{L}, \bigr) =
\lambda\cdot\pi_L = \lambda\cdot\biggl(\frac{\lambda}{1-p}
\biggr)^L\cdot\biggl(\frac{1-(\lambda/(1-p))}{1- (\lambda/(1-p))^{L+1}} \biggr). \label{eqrdthr}
\end{equation}
Define
%
%
\begin{equation}
L(x,\lambda)=\min\bigl\{L \in\mathbb{Z}_{+}\dvtx  r_d
\bigl(p,\lambda,\pi_{\mathrm{th}}^{L}, \bigr) \leq x \bigr\},
\label{eqLlamp}
\end{equation}
that is, $L(x,\lambda)$ is the smallest $L$ for which $\pi_{\mathrm{th}}^{L}$
remains feasible, given a deletion rate constraint of $x$. Using
equations~(\ref{eqrdthr}) and (\ref{eqLlamp}) to solve for
$L(p,\lambda)$, we obtain, after some algebra,
%
%
\begin{equation}\label{eqLlamp2}
L (p,\lambda) = \biggl\lceil\log_{\lambda/(1-p)}\frac{p}{1-\lambda} \biggr\rceil
\sim\log_{1/(1-\tilde{p})}\frac{1}{1-\lambda}\qquad\mbox{as } \lambda\to1
\end{equation}
and, by combining equation~(\ref{eqLlamp2}) and equation~(\ref
{eqCthresh}) with $L=L (p,\lambda)$, we have
%
%
\begin{equation}
C\bigl(p,\lambda, \pi_{\mathrm{th}}^{L(p,\lambda)}\bigr) \sim L(p,\lambda) \sim
\log_{1/(1-p)}\frac{1}{1-\lambda}\qquad\mbox{as } \lambda\to1. \label
{eqCthropt}
\end{equation}

By equations~(\ref{eqCmonL}) and (\ref{eqLlamp}), we know that
$\pi_{\mathrm{th}}^{L(p,\lambda)}$ achieves the minimum average queue length
among all feasible threshold policies. By equation~(\ref{eqCCC1}), we
must have that
%
%
\begin{equation}
C \bigl(p,\lambda,\pi_{\mathrm{th}}^{L (p,\lambda) -1} \bigr)\leq C\bigl
(p,\lambda,
\pi^*_p\bigr) \leq C \bigl(p,\lambda,\pi_{\mathrm{th}}^{L(p,\lambda)}
\bigr). \label{eqCCC2}
\end{equation}

Since Lemma~\ref{lemquasithresh} only applies when $B<\infty$,
equation~(\ref{eqCCC2}) holds whenever the buffer size, $B$, is
greater than $L(p,\lambda)$, but finite. We next extend equation~(\ref
{eqCCC2}) to the case of $B=\infty$. Denote by $\nu_p^*$ a
stationary optimal policy, when $B=\infty$ and the constraint on
average deletion rate is equal to $p$ [equation~(\ref{eqavgcost})].
The upper bound on $C(p,\lambda,\pi^*_p)$ in equation~(\ref
{eqCCC2}) automatically holds for $C(p,\lambda,\nu^*_p)$, since
$C(p,\lambda,\pi_{\mathrm{th}}^{L(p,\lambda)})$ is still feasible when
$B=\infty$. It remains to show a lower bound of the form
%
%
\begin{equation}
C\bigl(p,\lambda,\nu^*_p\bigr)\geq C \bigl(p,\lambda,
\pi_{\mathrm{th}}^{L
(p,\lambda) -2} \bigr), \label{eqClower}
\end{equation}
when\vspace*{1pt} $B=\infty$, which, together with the upper bound, will have
implied that the scaling of $ C(p,\lambda, \pi_{\mathrm{th}}^{L(p,\lambda)})$
[equation (\ref{eqCthropt})] carries over to $\nu^*_p$,
%
%
\begin{equation}
\qquad C \bigl(p,\lambda, \nu^*_p \bigr) \sim C\bigl(p,\lambda, \pi
_{\mathrm{th}}^{L(p,\lambda)}\bigr) \sim\log_{1/(1-p)}\frac{1}{1-\lambda}\qquad\mbox{as } \lambda\to1,
\end{equation}
thus proving Theorem~\ref{teoonline}.

To show equation~(\ref{eqClower}), we will use a straightforward
truncation argument that relates the performance of an optimal policy
under $B=\infty$ to the case of $B<\infty$. Denote by $ \{
b^*_q \}$ the deletion probabilities of a stationary optimal
policy, $\nu_p^*$, and by $ \{b^*_q(B') \}$ the deletion
probabilities for a truncated version, $\nu_p^*(B')$, with
\[
b^*_q\bigl(B'\bigr)=\mathbb{I} \bigl(q\leq
B' \bigr)\cdot b_q^*
\]
for all $q\geq0$. Since $\nu_p^*$ is optimal and yields the minimum
average queue length, it is without loss of generality to assume that
the Markov process for $Q(t)$ induced by $\nu_p^*$ is positive
recurrent. Denoting by $ \{\mu_i^* \}$ and $ \{\mu
^*_i(B') \}$ the steady-state probability of queue length being
equal to $i$ under $\nu_p^*$ and $\nu_p^*(B')$, respectively, it
follows from the positive recurrence of $Q(t)$ under $\nu_p$ and some
algebra, that
%
%
\begin{equation}
\lim_{B'\to\infty} \mu^{*}_i
\bigl(B'\bigr) = \mu^*_i \label{eqptwmucon}
\end{equation}
for all $i\in\mathbb{Z}_{+}$ and
%
%
\begin{equation}
\lim_{B'\to\infty} C \bigl(p,\lambda,\nu_p^*
\bigl(B'\bigr) \bigr) = C \bigl(p,\lambda,\nu_p^*
\bigr).
\end{equation}
By equation (\ref{eqptwmucon}) and the fact that $b_i^*(B')=b^*_i$
for all $0\leq i \leq B'$, we have that\footnote{Note that, in
general, $r_d (p,\lambda,\nu_p^*(B') )$ could be greater
than $p$, for any finite $B'$.}
%
%
\begin{eqnarray}
\lim_{B'\to\infty}r_d \bigl(p,\lambda,
\nu_p^*\bigl(B'\bigr) \bigr) &=& \lim
_{B'\to\infty}\lambda\sum_{i=0}^\infty
\mu^*_i\bigl(B'\bigr)\cdot\bigl(1-b_i^*
\bigl(B'\bigr) \bigr)
\nonumber
\\
&=& r_d \bigl(p,\lambda,
\nu_p^* \bigr)
\\
&\leq& p.\nonumber
\end{eqnarray}

It is not difficult to verify, from the definition of $L(p,\lambda)$
[equation~(\ref{eqLlamp})], that
\[
\lim_{\delta\to0} L(p+\delta,\lambda) \geq L(p,\lambda)-1,
\]
for all $p,\lambda$. For all $\delta>0$, choose $B'$ to be
sufficiently large, so that
%
%
\begin{eqnarray}
\label{eqCtoCB} C \bigl(p,\lambda,\nu_p^*\bigl(B'
\bigr) \bigr) &\leq& C \bigl(p,\lambda,\nu_p^* \bigr)+\delta,
\\
L \bigl(\lambda,r_d \bigl(p,\lambda,\nu_p^*
\bigl(B'\bigr) \bigr) \bigr) &\geq& L(p,\lambda)-1. \label{eqrdeps}
\end{eqnarray}

Let $p'=r_d (p,\lambda,\nu_p^*(B') )$. Since $b_i^*(B')=0$
for all $i\geq B'+1$, by equation~(\ref{eqrdeps}) we have
%
%
\begin{equation}
C \bigl(p,\lambda,\nu_p^*\bigl(B'\bigr) \bigr) \geq C
\bigl(p,\lambda,\pi^*_{p'} \bigr), \label{eqBinftoBfin}
\end{equation}
where $\pi^*_p$ is the optimal stationary policy given in Lemma~\ref{lemquasithresh} under any the finite buffer size $B>B'$. We have
%
%
\begin{eqnarray}\label{eqClower2}
&& C \bigl(p,\lambda,\nu_p^* \bigr)+\delta
\nonumber
\\
&&\qquad \stackrel{\mathrm{(a)}} {\geq}  C \bigl(p,\lambda,\nu_p^*
\bigl(B'\bigr) \bigr)
\nonumber
\\
&&\qquad\stackrel{\mathrm{(b)}} {\geq}  C \bigl(p,\lambda,\pi^*_{p'} \bigr)
\\
&&\qquad \stackrel{\mathrm{(c)}} {\geq}  C \bigl(p,\lambda,\pi_{\mathrm{th}}^{L(p',\lambda
)-1}\bigr)\nonumber
\\
&&\qquad \stackrel{\mathrm{(d)}} {\geq}  C \bigl(p,\lambda,\pi_{\mathrm{th}}^{L(p,\lambda
)-2}\bigr), \nonumber
\end{eqnarray}
where the inequalities (a) through (d) follow from equations~(\ref
{eqCtoCB}), (\ref{eqBinftoBfin}), (\ref{eqCCC2}) and (\ref
{eqrdeps}), respectively. Since equation~(\ref{eqClower2}) holds
for all $\delta>0$, we have proven equation~(\ref{eqClower}). This
completes the proof of Theorem~\ref{teoonline}.

\section{\texorpdfstring{Optimal offline policies.}{Optimal offline policies}}\label{secoffline}

We prove Theorem~\ref{teooffline} in this section, which is completed
in two parts. In the first part (Section~\ref{secperformoffline}),
we give a full characterization of the sample path resulted by applying
$\pi_{\mathrm{NOB}}$ (Proposition~\ref{propQRW}), which turns out to be a
\emph{recurrent} random walk. This allows us to obtain the
steady-state distribution of the queue length under $\pi_{\mathrm{NOB}}$ in
closed-form. From this, the expected queue length, which is equal to
the time-average queue length, $C (p,\lambda, \pi_{\mathrm{NOB}}
)$, can be easily derived and is shown to be $\frac{1-p }{\lambda
-(1-p)}$. Several side results we obtain along this path will also be
used in subsequent sections.

The second part of the proof (Section~\ref{secoptmalityofflineproof}) focuses on showing the heavy-traffic
optimality of $\pi_{\mathrm{NOB}}$ among the class of all feasible offline
policies, namely, that $\lim_{\lambda\to1} C (p,\lambda, \pi
_{\mathrm{NOB}} )=\lim_{\lambda\to1} C^*_{\Pi_\infty} (p,\lambda
)$, which, together with the first part, proves Theorem~\ref{teooffline} (Section~\ref{secpfthmoffline}). The optimality
result is proved using a sample-path-based analysis, by relating the
resulting queue length sample path of $\pi_{\mathrm{NOB}}$ to that of a greedy
deletion rule, which has an optimal deletion performance over a \emph
{finite} time horizon, $ \{1,\ldots, N \}$, given any
initial sample path. We then show that the discrepancy between $\pi
_{\mathrm{NOB}}$ and the greedy policy, in terms of the resulting time-average
queue length after deletion, diminishes almost surely as $N\to\infty$
and $\lambda\to1$ (with the two limits taken in this order). This
establishes the heavy-traffic optimality of $\pi_{\mathrm{NOB}}$.

\subsection{\texorpdfstring{Additional notation.}{Additional notation}}\label{sec6.1} Define $\widetilde{Q}$ as the resulting
queue length process after applying $\pi_{\mathrm{NOB}}$
\[
\widetilde{Q}=D \bigl(Q^0,M^\Psi\bigr) \label{eqqtil}
\]
and $Q$ as the shifted version of $\widetilde{Q}$, so that $Q$ starts
from the
first deletion in~$\widetilde{Q}$,\looseness=-1
%
%
\begin{equation}
Q[n] = \widetilde{Q}\bigl[n +m^\Psi_1\bigr],\qquad n \in
\mathbb{Z}_{+}.\vadjust{\goodbreak}
\end{equation}\looseness=0

We say that $B = \{l,\ldots,u \} \subset\mathbb{N}$ is a
\textit{busy period} of $Q$ if
%
%
\begin{equation}
\qquad Q[l-1]=Q[u]=0\quad\mbox{and}\quad Q[n]>0\qquad\mbox{for all }n\in\{l,\ldots,u-1 \}.
\label{eqbusydef}
\end{equation}
We may write $B_j = \{l_j,\ldots, u_j\}$ to mean the $j$th busy period
of $Q$. An example of a busy period is illustrated in Figure~\ref{figwater2}.

Finally, we will refer to the set of slots between two adjacent
deletions in $Q$ (note the offset of $m_1$),
%
%
\begin{equation}
E_i = \bigl\lbrace m^\Psi_i-m^\Psi_1,m^\Psi_i+1-m^\Psi
_1,\ldots,m^\Psi_{i+1}-1-m^\Psi_1
\bigr\rbrace \label{eqEjdef}
\end{equation}
as the $i$th \textit{deletion epoch}.

\subsection{\texorpdfstring{Performance of the no-job-left-behind policy.}{Performance of the no-job-left-behind policy}}\label{secperformoffline}

For simplicity of notation, throughout this section, we will denote by
$M= \{m_i\dvtx i \in\mathbb{N} \}$ the deletion sequence
generated by applying $\pi_{\mathrm{NOB}}$ to $Q^0$, when there is no ambiguity
(as opposed to using $M^\Psi$ and $m^\Psi_i$). The following lemma
summarizes some important properties of $Q$ which will be used repeatedly.

%
\begin{lem}
\label{lemqmi0}
Suppose $1>\lambda>1-p>0$. The following hold with probability one:
\begin{longlist}[(2)]
\item[(1)] For all $n\in\mathbb{N}$, we have $Q[n] = Q^0[n+m_1] -
I(M,n+m_1)$.
\item[(2)] For all $i\in\mathbb{N}$, we have $n=m_i-m_1$, if and
only if
%
%
\begin{equation}
Q[n]=Q[n-1]=0
\end{equation}
with the convention that $Q[-1]= 0$. In other words, the appearance of
two consecutive zeros in $Q$ is equivalent to having a deletion on the
second zero.
\item[(3)] $Q[n]\in\mathbb{Z}_{+}$ for all $n\in\mathbb{Z}_{+}$.
\end{longlist}
\end{lem}

\begin{pf}
See Appendix~\ref{applemqmi0}
\end{pf}

The next proposition is the main result of this subsection. It
specifies the probability law that governs the evolution of $Q$.

%
\begin{prop}
\label{propQRW}
$ \{Q[n]\dvtx n\in\mathbb{Z}_{+} \}$ is a random walk on
$\mathbb{Z}_{+}$, with $Q[0]=0$, and, for all $n\in\mathbb{N}$ and
$x_1,x_2\in\mathbb{Z}_{+}$,
\[
\mathbb{P} \bigl(Q[n+1]=x_2 \mid Q[n]=x_2 \bigr) =
\cases{ \displaystyle\frac{1-p}{\lambda+1-p}, &\quad$x_2-x_1=1$,
\vspace*{5pt}\cr
\displaystyle\frac{\lambda}{\lambda+1-p}, &\quad$x_2-x_1=-1$,
\vspace*{5pt}\cr
0, &\quad otherwise,}
\]
if $x_1>0$ and
\[
\mathbb{P} \bigl(Q[n+1]=x_2 \mid Q[n]=x_1 \bigr) =
\cases{ \displaystyle\frac{1-p}{\lambda+1-p}, &\quad$x_2-x_1=1$,
\vspace*{5pt}\cr
\displaystyle\frac{\lambda}{\lambda+1-p}, &\quad$x_2-x_1=0$,
\vspace*{5pt}\cr
0, &\quad otherwise,}
\]
if $x_1=0$.
\end{prop}

\begin{pf}
For a sequence $ \{X[n]\dvtx n\in\mathbb{N} \}$ and $s,t\in
\mathbb{N}$, $s\leq t$, we will use the shorthand
\[
X_{s}^t= \bigl\{X[s],\ldots, X[t] \bigr\}.
\]
Fix $n\in N$, and a sequence $ (q_1,\ldots,q_n )\subset
\mathbb{Z}_{+}^n$. We have
%
%
\begin{eqnarray}\label{eqQcond0}
&&\mathbb{P} \bigl(Q[n]=q[n] | Q_1^{n-1}=q_1^{n-1}
\bigr)\nonumber
\\
&&\qquad = \sum_{k=1}^{n} \mathop{\sum_{t_1,\ldots,t_k,}}_{t_k\leq
n-1+t_1} \mathbb{P} \bigl(Q[n]=q[n] | Q_1^{n-1}=q_1^{n-1},
m_{1}^k = t_1^k,
m_{k+1}\geq n +t_1 \bigr)\hspace*{-20pt}
\\
&&\hspace*{88pt}{}\times \mathbb{P} \bigl(m_{1}^k = t_1^k,
m_{k+1}\geq n +t_1 | Q_1^{n-1}=q_1^{n-1}
\bigr). \nonumber
\end{eqnarray}

Restricting to the values of $t_i$'s and $q[i]$'s under which the
summand is nonzero, the first factor in the summand can be written as
%
%
\begin{eqnarray}\label{eqQcond1}
\qquad&& \mathbb{P} \bigl(Q[n]=q[n] | Q_1^{n-1}=q_1^{n-1},
m_{1}^k = t_1^k,
m_{k+1}\geq n+t_1 \bigr)
\nonumber
\\
&&\qquad = \mathbb{P} \bigl(\widetilde{Q}[n+m_1]=q[n] | {\widetilde{Q}
}_{m_1+1}^{m_1+n-1} =q_1^{n-1},
m_{1}^k = t_1^k,
m_{k+1}\geq n+t_1 \bigr)\nonumber
\\
&&\qquad \stackrel{\mathrm{(a)}} {=}  \mathbb{P} \Bigl(Q^0[n+t_1]=q[n]+k
| Q^0[s+t_1]=q[s]+I \bigl( \{t_i
\}_{i=1}^k,s+t_1 \bigr),
\nonumber\\[-8pt]\\[-8pt]
&&\hspace*{153pt}\forall1\leq s \leq
n-1 \mbox{ and } \min_{r\geq n+t_1} Q^0[r] \geq k \Bigr)
\nonumber
\\
&&\qquad \stackrel{\mathrm{(b)}} {=}  \mathbb{P} \Bigl(Q^0[n+t_1]=q[n]+k
| Q^0[n-1+t_1]=q[n-1]+k\nonumber
\\
&&\hspace*{190.5pt} \mbox{and } \min_{r\geq n+t_1} Q^0[r] \geq k \Bigr), \nonumber
\end{eqnarray}
where $\widetilde{Q}$ was defined in equation~(\ref{eqqtil}). Step~(a)
follows from Lemma~\ref{lemqmi0} and the fact that $t_k\leq n-1+t_1$,
and (b) from the Markov property of $Q^0$ and the fact that the
events $ \{\min_{r\geq n+t_1} Q^0[r] \geq k \}$, $ \{
Q^0[n+t_1]=q[n]+k \}$ and their intersection, depend only on the
values of $ \{Q^0[s]\dvtx s\geq n+t_1 \}$, and are hence
independent of $ \{Q^0[s]\dvtx 1\leq s \leq n-2+t_1 \}$
conditional on the value of $Q^0[t_1+n-1]$.

Since the process $Q$ lives in $\mathbb{Z}_{+}$ (Lemma~\ref{lemqmi0}), it suffices to consider the case of $q[n]=q[n-1]+1$, and
show that
%
%
\begin{eqnarray}
&& \mathbb{P} \Bigl(Q^0[n+t_1]=q[n-1]+1+k |
Q^0[n-1+t_1]=q[n-1]+k \nonumber
\\
&&\hspace*{192pt} \mbox{and }\min _{r\geq n+t_1} Q^0[r] \geq k \Bigr)
\\
&&\qquad = \frac{1-p}{\lambda+1-p}\nonumber
\end{eqnarray}
for all $q[n-1]\in\mathbb{Z}_{+}$. Since $Q[m_i-m_1]=Q[m_i-1-m_1]=0$
for all $i$ (Lemma~\ref{lemqmi0}), the fact that $q[n]=q[n-1]+1>0$
implies that
%
%
\begin{equation}
n < m_{k+1}-1+m_1. \label{eqn1}
\end{equation}
Moreover, since $Q^0[m_{k+1}-1]=k$ and $n< m_{k+1}-1+m_1$, we have that
%
%
\begin{equation}
q[n]>0\qquad\mbox{implies } Q^0[t]=k\qquad\mbox{for some } t\geq n+1+m_1. \label{eqAcondition}
\end{equation}
We consider two cases, depending on the value of $q[n-1]$.
\begin{longlist}
\item[{\textit{Case} 1: $q[n-1]>0$.}] Using the same argument that led to
equation~(\ref{eqAcondition}), we have that
%
%
\begin{equation}
q[n-1]>0\qquad\mbox{implies } Q^0[t]=k\qquad\mbox{for some } t\geq n
+m_1. \label{eqn2}
\end{equation}
It is important to note that, despite the similarity of their
conclusions, equations~(\ref{eqAcondition}) and (\ref{eqn2}) are
different in their assumptions (i.e., $q[n]$ versus $q[n-1]$). We have
%
%
\begin{eqnarray}\label{eqQcond2}
\qquad &&\mathbb{P} \Bigl(Q^0[n+t_1]=q[n-1]+1+k |
Q^0[n-1+t_1]=q[n-1]+k\nonumber
\\
&&\hspace*{193pt}\mbox{and } \min
_{r\geq n+t_1} Q^0[r] \geq k \Bigr)\nonumber
\\
&&\qquad \stackrel{\mathrm{(a)}} {=}\mathbb{P} \Bigl(Q^0[n+t_1]=q[n-1]+1+k
| Q^0[n-1+t_1]=q[n-1]+k
\nonumber\\[-8pt]\\[-8pt]
&&\hspace*{227pt}\mbox{and } \min _{r\geq n+t_1} Q^0[r] = k \Bigr)\nonumber
\\
&&\qquad \stackrel{\mathrm{(b)}} {=}\mathbb{P} \Bigl(Q^0[2]=q[n-1]+1 |
Q^0[1]=q[n-1]\mbox{ and } \min_{r\geq2}Q^0[r] = 0 \Bigr)\nonumber
\\
&&\qquad \stackrel{\mathrm{(c)}} {=}  \frac{1-p}{\lambda+1-p}, \nonumber
\end{eqnarray}
where (a) follows from equation~(\ref{eqn2}), (b) from the
stationary and space-homogeneity of the Markov chain $Q^0$ and (c)
from the following well-known property of a transient random walk
conditional to returning to zero:
\end{longlist}
%
\begin{lem}
\label{lemdualRW}
Let $ \{X[n]\dvtx n\in\mathbb{N} \}$ be a random walk on
$\mathbb{Z}_{+}$, such that for all $x_1,x_2\in\mathbb{Z}_{+}$ and
$n\in\mathbb{N}$,
\[
\mathbb{P} \bigl(X[n+1]=x_2 \mid X[n]=x_2 \bigr) =
\cases{ q, &\quad$x_2-x_1=1$,
\vspace*{5pt}\cr
1-q, &
\quad$x_2-x_1=-1$,
\vspace*{5pt}\cr
0, &\quad otherwise,}
\]
if $x_1>0$ and
\[
\mathbb{P} \bigl(X[n+1]=x_2 \mid X[n]=x_1 \bigr) =
\cases{ q, &\quad$x_2-x_1=1$,
\vspace*{5pt}\cr
1-q, &
\quad$x_2-x_1=0$,
\vspace*{5pt}\cr
0, &\quad otherwise,}
\]
if $x_1=0$, where $q\in(\frac{1}{2},1 )$. Then for all
$x_1,x_2\in\mathbb{Z}_{+}$ and $n\in\mathbb{N}$,
\[
\mathbb{P} \Bigl(X[n+1]=x_2 | X[n]=x_1, \min
_{r \geq n+1} X[r] = 0 \Bigr) = \cases{ 1-q, &\quad$x_2-x_1=1$,
\vspace*{5pt}\cr
q, &\quad$x_2-x_1=-1$,
\vspace*{5pt}\cr
0, &\quad otherwise,}
\]
if $x_1>0$ and
\[
\mathbb{P} \Bigl(X[n+1]=x_2 | X[n]=x_1, \min
_{r \geq n+1} X[r] = 0 \Bigr) = \cases{ 1-q, &\quad$x_2-x_1=1$,
\vspace*{5pt}\cr
q, &\quad$x_2-x_1=0$,
\vspace*{5pt}\cr
0, &\quad otherwise,}
\]
if $x_1=0$. In other words, conditional on the eventual return to $0$
and before it happens, a transient random walk obeys the same
probability law as a random walk with the reversed one-step transition
probability.
\end{lem}

\begin{pf}
See Appendix~\ref{applemdualRW}.
\end{pf}
\begin{longlist}
\item[{\textit{Case} 2: $q[n-1]=0$}] We have
%
%
\begin{eqnarray}\label{eqQcond3}
\qquad &&\mathbb{P} \Bigl(Q^0[n+t_1]=q[n-1]+1+k |
Q^0[n-1+t_1]=q[n-1]+k\nonumber
\\
&&\hspace*{190pt} \mbox{ and } \min _{r\geq n+t_1} Q^0[r] \geq k \Bigr)\nonumber
\\
&&\qquad \stackrel{\mathrm{(a)}} {=}\mathbb{P} \Bigl(Q^0[n+t_1]=1+k\nonumber
\mbox{ and } \min_{r> n+t_1} Q^0[r] = k |
Q^0[n-1+t_1]=k
\nonumber\\[-8pt]\\[-8pt]
&&\hspace*{220pt} \mbox{ and } \min_{r\geq n+t_1}
Q^0[r] \geq k \Bigr)\nonumber
\\
&&\qquad \stackrel{\mathrm{(b)}} {=} \mathbb{P} \Bigl(Q^0[2]=2 \mbox{ and } \min
_{r>2} Q^0[r] = 1 | Q^0[1]=1 \mbox{
and } \min_{r\geq2} Q^0[r] \geq1 \Bigr),
\nonumber
\\
&&\qquad \stackrel{\triangle} {=} x, \nonumber
\end{eqnarray}
where (a) follows from equation~(\ref{eqAcondition}) [note its
difference with equation~(\ref{eqn2})], and (b) from the
stationarity and space-homogeneity of $Q^0$, and the assumption that
$k\geq1$ [equation~(\ref{eqQcond0})].

Since equations~(\ref{eqQcond2}) and (\ref{eqQcond3}) hold for all
$x_1,k\in\mathbb{Z}_{+}$ and $n\geq m_1+1$, by equation~(\ref
{eqQcond0}), we have that
%
%
\begin{eqnarray}\label{eqtransit1}
&& \mathbb{P} \bigl(Q[n]=q[n] | Q_1^{n-1}=q_1^{n-1}
\bigr)
\nonumber\\[18pt]\\[-30pt]
&&\qquad = \cases{\displaystyle\frac{1-p}{\lambda+1-p}, &\quad$q[n]-q[n-1]=1$,
\vspace*{5pt}\cr
\displaystyle\frac{\lambda}{\lambda+1-p}, &\quad$q[n]-q[n-1]=-1$,
\vspace*{5pt}\cr
0, &\quad otherwise,}\nonumber
\end{eqnarray}
if $q[n-1]>0$ and
%
%
\begin{eqnarray}
\qquad\mathbb{P} \bigl(Q[n]=q[n] | Q_1^{n-1}=q_1^{n-1}
\bigr) = \cases{x, &\quad$q[n]-q[n-1]=1$,
\vspace*{5pt}\cr
1-x, &\quad$q[n]-q[n-1]=0$,
\vspace*{5pt}\cr
0, &
\quad otherwise,}\label{eqtransit2}
\end{eqnarray}
if $q[n-1]=0$, where $x$ represents the value of the probability in
equation~(\ref{eqQcond3}). Clearly, $Q[0]=Q^0[m_1]=0$. We next show
that $x$ is indeed equal to $\frac{1-p}{\lambda+1-p}$, which will
have proven Proposition~\ref{propQRW}.

One can in principle obtain the value of $x$ by directly computing the
probability in line (b) of equation~(\ref{eqQcond3}), which can be
quite difficult to do. Instead, we will use an indirect approach that
turns out to be computationally much simpler: we will relate $x$ to the
rate of deletion of $\pi_{\mathrm{NOB}}$ using renewal theory, and then solve
for $x$. As a by-product of this approach, we will also get a better
understanding of an important regenerative structure of $\pi_{\mathrm{NOB}}$
[equation~(\ref{eqlimfrac1})], which will be useful for the analysis
in subsequent sections.

By equations~(\ref{eqtransit1}) and (\ref{eqtransit2}), $Q$ is a
positive recurrent Markov chain, and $Q[n]$ converges to a well-defined
steady-state distribution, $Q[\infty]$, as $n\to\infty$. Letting
$\pi_i = \mathbb{P} (Q[\infty]=i )$, it is easy to verify
via the balancing equations that
%
%
\begin{equation}
\pi_i = \pi_0\frac{x(\lambda+1-p)}{\lambda}\cdot\biggl(
\frac
{1-p}{\lambda} \biggr)^{i-1}\qquad\forall i\geq1
\end{equation}
and since $\sum_{i\geq0}\pi_i=1$, we obtain
%
%
\begin{equation}
\pi_0 = \frac{1}{1+x\cdot (\lambda+1-p)/(\lambda-(1-p))}.
\end{equation}
Since the chain $Q$ is also irreducible, the limiting fraction of time
that $Q$ spends in state 0 is therefore equal to $\pi_0$,
%
%
\begin{equation}
\qquad \lim_{n\to\infty} \frac{1}{n}\sum
_{t=1}^n \mathbb{I} \bigl(Q[t]=0 \bigr) =
\pi_0= \frac{1}{1+x\cdot(\lambda+1-p)/(\lambda-(1-p))}. \label{eqlimzero}
\end{equation}

Next, we would like to know many of these visits to state 0 correspond
to a deletion. Recall the notion of a busy period and deletion epoch,
defined in equations~(\ref{eqbusydef}) and (\ref{eqEjdef}),
respectively. By Lemma~\ref{lemqmi0}, $n$ corresponds to a deletion
if any only if $Q[n]=Q[n-1]=0$. Consider a deletion in slot $m_i$. If
$Q[m_i+1]=0$, then $m_i+1$ also corresponds to a deletion, that is,
$m_i+1 = m_{i+1}$. If instead $Q[m_i+1]=1$, which happens with
probability $x$, the fact that $Q[m_{i+1}-1]=0$ implies that there
exists at least one busy period, $\{l,\ldots,u\}$, between $m_i$ and
$m_{i+1}$, with $l=m_i$ and $u \leq m_{i+1}-1$. At the end of this
period, a new busy period starts with probability $x$ and so on. In
summary, a deletion epoch $E_i$ consists of the slot $m_i-m_1$, plus
$N_i$ busy periods, where the $N_i$ are i.i.d., with\footnote{$\operatorname{Geo}(p)$ denotes a geometric random variable with mean $\frac{1}{p}$.}
%
%
\begin{equation}
N_1 \stackrel{d} {=} \operatorname{Geo}(1-x)-1
\end{equation}
and hence
%
%
\begin{equation}\label{eqEidecomp}
|E_i| = 1+\sum_{j=1}^{N_i}
B_{i,j},
\end{equation}
where $ \{B_{i,j}\dvtx i,j\in\mathbb{N} \}$ are i.i.d. random
variables, and $B_{i,j}$ corresponds to the length of the $j$th busy
period in the $i$th epoch.

Define $W[t] = (Q[t],Q[t+1] )$, $t\in\mathbb{Z}_{+}$.
Since $Q$ is Markov, $W[t]$ is also a Markov chain, taking values in
$\mathbb{Z}_{+}^2$. Since a deletion occurs in slot $t$ if and only if
$Q[t]=Q[t-1]=0$ (Lemma~\ref{lemqmi0}), $|E_i|$ corresponds to
excursion times between two adjacent visits of $W$ to the state
$(0,0)$, and hence are i.i.d. Using the elementary renewal theorem, we have
%
%
\begin{equation}
\lim_{n\to\infty} \frac{1}{n} I (M,n ) = \frac
{1}{\mathbb{E}(|E_1|)}\qquad\mbox{a.s.} \label{eqlimfrac1}
\end{equation}
and by viewing each visit of $W$ to $(0,0)$ as a renewal event and
using the fact that exactly one deletion occurs within a deletion
epoch. Denoting by $R_i$ the number of visits to the state 0 within
$E_i$, we have that $R_i=1+N_i$. Treating $R_i$ as the reward
associated with the renewal interval $E_i$, we have, by the
time-average of a renewal reward process (cf. Theorem 6, Chapter~3,
\cite{Gal96}), that
%
%
\begin{equation}
\lim_{n\to\infty} \frac{1}{n} \sum
_{t=1}^n \mathbb{I} \bigl(Q[t]=0 \bigr) =
\frac{\mathbb{E} (R_1 )}{\mathbb
{E} (|E_1| )} = \frac{\mathbb{E} (N_1
)+1}{\mathbb{E} (|E_1| )}\qquad\mbox{a.s.} \label{eqlimfrac2}
\end{equation}
by treating each visit of $Q$ to $(0,0)$ as a renewal event. From
equations (\ref{eqlimfrac1})~and~(\ref{eqlimfrac2}), we have
%
%
\begin{equation}
\frac{\lim_{n\to\infty} (1/n)I (M,n )}{\lim_{n\to\infty} (1/n) \sum_{t=1}^n \mathbb{I}
(Q[t]=0 ) } = \frac{1}{\mathbb{E}(N_1)} = 1-x. \label{eqlimfrac}
\end{equation}
Combining equations~(\ref{eqIMnlim}), (\ref{eqlimzero}) and (\ref
{eqlimfrac}), and the fact that $\mathbb{E}(N_1)=\mathbb{E}(\operatorname{Geo}(1-x))-1=\frac{1}{1-x}-1$, we have
%
%
\begin{equation}
\frac{\lambda-(1-p)}{\lambda+1-p} \cdot\biggl[1+x\cdot\frac
{\lambda+1-p}{\lambda-(1-p)} \biggr]= 1-x,
\end{equation}
which yields
%
%
\begin{equation}
x=\frac{1-p}{\lambda+1-p}.
\end{equation}

This completes the proof of Proposition~\ref{propQRW}.\quad\qed
\end{longlist}\noqed
\end{pf}

We summarize some of the key consequences of Proposition~\ref{propQRW} below, most of which are easy to derive using renewal theory
and well-known properties of positive-recurrent random walks.

%
\begin{prop}
\label{proppnobperf}
Suppose that $1>\lambda>1-p>0$, and denote by $Q[\infty]$ the
steady-state distribution of $Q$.
\begin{longlist}[(2)]
\item[(1)] For all $i\in\mathbb{Z}_{+}$,
%
%
\begin{equation}
\mathbb{P} \bigl(Q[\infty]=i \bigr) = \biggl(1-\frac{1-p}{\lambda
} \biggr)\cdot
\biggl(\frac{1-p}{\lambda} \biggr)^i.
\end{equation}
\item[(2)] Almost surely, we have that
%
%
\begin{equation}
\lim_{n\to\infty}\frac{1}{n}\sum_{i=1}^n
Q[i] = \mathbb{E} \bigl(Q [\infty] \bigr) = \frac{1-p}{\lambda-
(1-p )}. \label{eqQavg}
\end{equation}
\item[(3)] Let $E_i = \{m^\Psi_{i},m^\Psi_{i}+1,\ldots,m^\Psi_{i+1}-1,m^\Psi_{i+1} \}$. Then the $|E_i|$ are
i.i.d., with
%
%
\begin{equation}
\mathbb{E} \bigl(|E_1| \bigr) = \frac{1}{\lim_{n\to\infty}(1/n) I (M^\Psi,n )}= \frac{\lambda+1-p}{\lambda-(1-p)}
\end{equation}
and there exists $a,b>0$ such that for all $x\in\mathbb{R}_+$
%
%
\begin{equation}
\mathbb{P} \bigl(|E_1| \geq x \bigr) \leq a\cdot\exp(- b\cdot x ).
\label{eqEiExp}
\end{equation}
\item[(4)] Almost surely, we have that
%
%
\begin{equation}
m^\Psi_{i} \sim\frac{1}{\mathbb{E} (|E_1| )}\cdot i =
\frac{\lambda-(1-p)}{\lambda+1-p}\cdot i \label{eqmilim}
\end{equation}
as $i \to\infty$.
\end{longlist}
\end{prop}
\begin{pf} Claim $1$ follows from the well-known steady-state
distribution of a random walk, or equivalently, the fact that $Q[\infty
]$ has the same distribution as the steady-state number of jobs in an
$M/M/1$ queue with traffic intensity $\rho= \frac{1-p}{\lambda}$.
For Claim $2$, since $Q$ is an irreducible Markov chain that is
positive recurrent, it follows that its time-average coincides with
$\mathbb{E} (Q[\infty] )$ almost surely.

The fact that $E_i$'s are i.i.d. was shown in the discussion preceding
equation~(\ref{eqlimfrac1}) in the proof of Proposition~\ref{propQRW}. The value of $\mathbb{E} (|E_1| )$ follows by
combining equations~(\ref{eqIMnlim}) and (\ref{eqlimfrac1}).

Let $B_{i,j}$ be the length of the $j$th busy period [defined in
equation~(\ref{eqbusydef})] in $E_i$. By definition, $B_{1,1}$ is
distributed as the time till the random walk $Q$ reaches state $0$,
starting from state $1$. We have
\[
\mathbb{P} (B_{1,1}\geq x ) \leq\mathbb{P} \Biggl(\sum
_{j=1}^{ \lfloor x \rfloor} X_j \leq-1 \Biggr),
\]
where the $X_j$'s are i.i.d., with $\mathbb{P} (X_1=1
)=\frac{1-p}{\lambda+1-p}$ and $\mathbb{P} (X_1=-1
)=\frac{\lambda}{\lambda+1-p}$, which, by the Chernoff bound,
implies an exponential tail bound for\break  $\mathbb{P} (B_{1,1}\geq x
)$, and in particular,
%
%
\begin{equation}\label{eqGBlim}
\lim_{\theta\downarrow0} G_{B_{1,1}}(\theta) = 1. \vadjust{\goodbreak}
\end{equation}
By equation~(\ref{eqEidecomp}), the moment generating function for
$|E_1|$ is given by
%
%
\begin{eqnarray}
G_{|E_1|}(\varepsilon) &= & \mathbb{E} \bigl(\exp\bigl(\varepsilon\cdot
|E_1| \bigr) \bigr)
\nonumber
\\
&=& \mathbb{E} \Biggl(\exp\Biggl(\varepsilon\cdot\Biggl(1+\sum
_{j=1}^{N_1}B_{1,j} \Biggr) \Biggr) \Biggr)
\nonumber\\[-8pt]\\[-8pt]
&\stackrel{\mathrm{(a)}} {=} & \mathbb{E} \bigl(e^\varepsilon\bigr)\cdot\mathbb{E}
\bigl(\exp\bigl({N_1}\cdot G_{B_{1,1}}(\varepsilon) \bigr) \bigr)
\nonumber
\\
&= & \mathbb{E} \bigl(e^\varepsilon\bigr)\cdot G_{N_1} \bigl(\ln
\bigl(G_{B_{1,1}}(\varepsilon) \bigr) \bigr),\nonumber
\end{eqnarray}
where (a) follows from the fact that $ \{N_1 \}\cup
\{B_{1,j}\dvtx j\in\mathbb{N} \}$ are mutually independent, and
$G_{N_1}(x)= \mathbb{E} (\exp(x\cdot N_1 ) )$.
Since $N_1 \stackrel{d}{=} \operatorname{Geo}(1-x)-1$,\break  $\lim_{x\downarrow
0}G_{N_1}(x)=1$, and by equation~(\ref{eqGBlim}), we have that $\lim
_{\varepsilon\downarrow0}G_{|E_1|}(\varepsilon)=1$, which implies
equation~(\ref{eqEiExp}).

Finally, equation~(\ref{eqmilim}) follows from the third claim and
the elementary renewal theorem.
\end{pf}

\subsection{\texorpdfstring{Optimality of the no-job-left-behind policy in heavy traffic.}{Optimality of the no-job-left-behind policy in heavy traffic}}\label{secoptmalityofflineproof}
This section is devoted to proving the optimality of $\pi_{\mathrm{NOB}}$ as
$\lambda\to1$, stated in the second claim of Theorem~\ref{teooffline}, which we isolate here in the form of the following proposition.
%
\begin{prop} \label{propnobopt} Fix $p\in(0,1)$. We have that
\[
\lim_{\lambda\rightarrow1}C (p,\lambda,\pi_{\mathrm{NOB}} )=\lim
_{\lambda\rightarrow1}C^*_{\Pi_\infty} (p,\lambda).
\]
\end{prop}
The proof is given at the end of this section, and we do so by showing
the following:
\begin{longlist}[(2)]
\item[(1)] Over a finite horizon $N$ and given a fixed number of
deletions to be made, a~greedy deletion rule is optimal in minimizing
the post-deletion area under $Q$ over $ \{1,\ldots, N \}$.
\item[(2)] Any point of deletion chosen by $\pi_{\mathrm{NOB}}$ will also be
chosen by the greedy policy, as $N\to\infty$.
\item[(3)] The fraction of points chosen by the greedy policy but not
by $\pi_{\mathrm{NOB}}$ diminishes as $\lambda\to1$, and hence the delay
produced by $\pi_{\mathrm{NOB}}$ is the best possible, as $\lambda\to1$.
\end{longlist}

Fix $N\in\mathbb{N}$. Let $S (Q,N )$ be the
partial sum $S (Q,N )=\sum_{n=1}^{N}Q [n ]$. For
any sample path $Q$,
denote by $\Delta(Q,n )$ the marginal decrease of area under
$Q$ over the horizon $ \{1,\ldots,N \}$ by applying a
deletion at
slot $n$, that is,
\[
\Delta_P (Q,N,n )=S (Q,N )-S \bigl(D_P (Q,n ),N
\bigr)
\]
and, analogously,
\[
\Delta\bigl(Q,N,M' \bigr)=S (Q,N )-S \bigl(D
\bigl(Q,M' \bigr),N \bigr),
\]
where $M'$ is a deletion sequence.

We next define the notion of a greedy deletion rule, which constructs a
deletion sequence by recursively adding the slot that leads to the
maximum marginal decrease in $S(Q,N)$.

%
\begin{defn}[(Greedy deletion rule)]\label{defgdRule}
Fix an initial sample path $Q^{0}$
and $K,N\in\mathbb{N}$.
The \textit{greedy deletion rule} is a mapping, $G
(Q^{0},N,K )$,
which outputs a finite deletion sequence $M^G= \{
m_{i}^{G}\dvtx 1\leq i\leq K \} $,
given by
\begin{eqnarray*}
m_{1}^{G} & \in& \arg\max_{m\in\Phi(Q^0,N )}\Delta
_P \bigl(Q^{0},N,m \bigr),
\\
m_{k}^{G} & \in& \arg\max_{m\in\Phi(Q^{k-1},N )}\Delta
_P \bigl(Q_{M^G}^{k-1},N,m \bigr),\qquad2\leq k
\leq K,
\end{eqnarray*}
where $\Phi(Q,N ) = \Phi(Q ) \cap\{
1,\ldots,N \}$ is the set of all locations in $Q$ in the first
$N$~slots that can be deleted, and $Q_{M^G}^{k}=D (Q^{0},
\{ m_{i}^{G}\dvtx 1\leq i\leq k \} )$. Note that we will allow
$m_{k}^{G}= \infty$,
if there is no more entry to delete [i.e., $\Phi(Q^{k-1}
)\cap\{1,\ldots,N \}=\varnothing$].
\end{defn}
We now state a key lemma that will be used in proving Theorem~\ref{teooffline}.
It shows that over a finite horizon and for a finite number of
deletions, the greedy deletion rule yields the maximum reduction in the
area under the sample path.

%
\begin{lem}[(Dominance of greedy policy)]\label{lemgddom}
Fix an initial
sample path $Q^{0}$, horizon $N \in\mathbb{N}$ and number of
deletions $K\in\mathbb{N}$.
Let $M'$ be any deletion sequence with $I(M',N)=K$. Then
\[
S \bigl(D \bigl(Q^{0},M' \bigr),N \bigr)\geq S \bigl(D
\bigl(Q^{0},M^G \bigr),N \bigr),
\]
where $M^G=G (Q^{0},N,K )$ is the deletion sequence
generated by the greedy policy.
\end{lem}
\begin{pf} By Lemma~\ref{lemlocindp}, it suffices to show
that, for any sample path $ \{Q[n]\in\mathbb{Z}_{+}\dvtx n\in\mathbb
{N} \}$ with $|Q[n+1]-Q[n]|=1$ if $Q[n]>0$ and $|Q[n+1]-Q[n]|\in
\{0,1 \}$ if $Q[n]=0$, we have
%
%
\begin{eqnarray}\label{eqgrrec}
&& S \bigl(D \bigl(Q,M' \bigr),N \bigr)
\nonumber\\[-8pt]\\[-8pt]
&&\qquad \geq\Delta_P \bigl(Q,N,m_{1}^{G} \bigr)+\mathop{\min_{\llvert\widetilde{M}\rrvert
=k-1,}}_{\widetilde{M}\subset\Phi(D (Q,m_{1}^{G}
),N )}S\bigl(D \bigl(Q_{M^G}^{1},\widetilde{M} \bigr),N
\bigr).\nonumber
\end{eqnarray}
By induction, this would imply that we should use the greedy
rule at every step of deletion up to $K$. The following lemma states
a simple monotonicity property. The proof is elementary,
and is omitted.
%
\begin{lem}[(Monotonicity in deletions)]\label{lemMonotonicity-in-Deletions}
Let $Q$ and $Q'$ be two sample paths such that
\[
Q [n ]\leq Q' [n ]\qquad\forall n\in\{ 1,\ldots,N \}.
\]
Then, for any $K\geq1$,
%
%
\begin{equation}\label{eqmonDel1}
\mathop{\min_{\llvert M\rrvert=K,}}_{M\subset\Phi(Q,N
)}S \bigl(D (Q,M ),N \bigr)\leq\mathop{\min_{\llvert M\rrvert=K,}}_{M\subset\Phi(Q',N )}S \bigl(D \bigl
(Q',M \bigr),N
\bigr)
\end{equation}
and, for any finite deletion sequence $M'\subset\Phi(Q,N )$,
%
%
\begin{equation}
\Delta\bigl(Q,N,M' \bigr)\geq\Delta\bigl(Q',N,M'
\bigr).\label
{eqmonDel2}
\end{equation}
\end{lem}

Recall the definition of a busy period in equation~(\ref{eqbusydef}).
Let $J(Q,N)$ be the total number of busy periods in $ \{Q[n]\dvtx 1\leq
n\leq N \}$, with the additional convention $Q[N+1]\stackrel
{\triangle}{=}0$ so that the last busy period always ends on $N$. Let
$B_j = \{l_j,\ldots,u_j \}$ be the $j$th busy period. It
can be verified that a deletion in location $n$ leads to a decrease in
the value of $S(Q,N)$ that is no more than the width of the busy period
to which $n$ belongs; cf. Figure~\ref{figwater2}. Therefore, by
definition, a greedy policy always seeks to delete in each step the
first arriving job during a longest busy period in the current sample
path, and hence
%
%
\begin{equation}
\Delta\bigl(Q,N, G(Q,N,1) \bigr) = \max_{1\leq j \leq J(Q,N)}\llvert
B_{j}\rrvert. \label{eqbubble1}
\end{equation}


Let
\[
\mathcal{J}^{*}(Q,N)= \arg\max_{1\leq j \leq J(Q,N)} \llvert
B_j\rrvert.
\]
We consider the following cases, depending on whether $M'$ chooses to
delete any job in the busy periods in $\mathcal{J}^{*}(Q,N)$.
\begin{longlist}
\item[{\textit{Case} 1: $M' \cap(\bigcup_{j\in\mathcal
{J}^{*}(Q,N)}B_{j} )\neq\varnothing$.}]
If $l_{j^*}\in M'$ for some $j^{*}\in\mathcal{J}^{*}$,
by equation~(\ref{eqbubble1}), we can set $m_{1}^{G}$ to $l_{j^*}$.
Since $m_{1}^{G}\in M'$ and the order of deletions does not impact the
final resulting delay (Lemma~\ref{lemlocindp}), we have that
equation~(\ref{eqgrrec}) holds, and we are done. Otherwise, choose
$m^{*}\in M'\cap B_{j^{*}}$ for some $j^{*}\in\mathcal{J}^{*}$, and
we have $m^*> l_{j^*}$. Let
\[
Q'= D_P \bigl(Q,m^{*} \bigr)\quad\mbox{and}\quad \widehat{Q} = D_P (Q,l_{j^*} ).
\]
Since $Q [n ]>0$, $\forall n\in\{l_{j^*}, \ldots,
u_{j^*}-1 \}$, we have $\widehat{Q} [n ]=Q
[n
]-1\leq Q'[n]$, $\forall n\in\{l_{j^*}, \ldots,
u_{j^*}-1 \}$ and $Q'[n]=Q[n]=\widehat{Q}[n]$, $\forall n \notin
\{l_{j^*}, \ldots, u_{j^*}-1 \}$, which implies that
%
%
\begin{equation}
\widehat{Q} [n ]\leq Q' [n ]\qquad\forall n\in\{1,\ldots,N
\}.\label{eqcase1dom}
\end{equation}
Equation~(\ref{eqgrrec}) holds by combining equation~(\ref{eqcase1dom})
and equation~(\ref{eqmonDel1}) in Lemma~\ref{lemMonotonicity-in-Deletions},
with $K=k-1$.
\end{longlist}
\begin{longlist}
\item[{\textit{Case} 2: $M'\cap(\bigcup_{j\in\mathcal{J}^{*}(Q,N)}B_{j} )= \varnothing$.}]
Let\vspace*{2pt} $m^{*}$ be any element in $M'$ and $Q'=D_P (Q,m^{*} )$.
Clearly, $Q [n ]\geq Q' [n ]$ for\vspace*{1pt} all $n\in
\{1,\ldots,N \}$, and by equation~(\ref{eqmonDel2}) in Lemma
\ref{lemMonotonicity-in-Deletions}, we have that\footnote{For
finite sets $A$ and $B$, $A\setminus B = \{a\in A\dvtx  a\notin
B \}$.}
%
%
\begin{equation}
\Delta\bigl(Q,N,M'\setminus\bigl\{ m^{*} \bigr\}
\bigr)\geq\Delta\bigl(D_P \bigl(Q,m^* \bigr),N,M'
\setminus\bigl\{ m^* \bigr\} \bigr). \label{eqcase2-diff-dom}
\end{equation}
Since \textbf{$M'\cap(\bigcup_{j\in\mathcal
{J}^{*}(Q,N)}B_{j} )=\varnothing$},
we have that
%
%
\begin{equation}
\qquad\Delta_P \bigl(D \bigl(Q,M'\setminus\bigl\{
m^{*} \bigr\} \bigr),N,m_{1}^{G} \bigr)=\max
_{1\leq j \leq J(Q,N)}\llvert B_{j}\rrvert>\Delta_P
\bigl(Q,N,m^{*} \bigr).\label{eqcase2-step-dom}
\end{equation}
Let $\widehat{M}=m_{1}^{G}\cup(M'\setminus\{
m^{*} \}
)$, and
we have that
\begin{eqnarray*}
&& S \bigl(D (Q,\widehat{M} ),N \bigr)
\\
&&\qquad =  S (Q,N )-\Delta\bigl(Q,N,M'\setminus\bigl\{
m^{*} \bigr\} \bigr)-\Delta_P \bigl(D
\bigl(Q,M'\setminus\bigl\{ m^{*} \bigr\}
\bigr),N,m_{1}^{G} \bigr)
\\
&&\qquad \stackrel{\mathrm{(a)}} {\leq}  S (Q,N )-\Delta\bigl(D_P \bigl(Q,m^*
\bigr),N,M'\setminus\bigl\{ m^{*} \bigr\} \bigr)
\\
&&\quad\qquad{} -\Delta_P \bigl(D \bigl(Q,M'\setminus\bigl\{
m^{*} \bigr\} \bigr),N,m_{1}^{G} \bigr)
\\
&&\qquad \stackrel{\mathrm{(b)}} {<}  S (Q,N )-\Delta\bigl(D_P \bigl(Q,m^*
\bigr),N,M'\setminus\bigl\{ m^{*} \bigr\} \bigr)-\Delta
_P \bigl(Q,N,m^{*} \bigr)
\\
&&\qquad = S \bigl(D \bigl(Q,M' \bigr),N \bigr),
\end{eqnarray*}
where (a) and (b) follow from equations~(\ref{eqcase2-diff-dom})
and (\ref{eqcase2-step-dom}), respectively, which shows that
equation~(\ref{eqgrrec}) holds (and in this case the inequality
there is strict).

Cases 1~and~2 together complete the proof of Lemma~\ref{lemgddom}.\quad\qed
\end{longlist}\noqed
\end{pf}

We are now ready to prove Proposition~\ref{propnobopt}.

\begin{pf*}{Proof of Proposition~\ref{propnobopt}}
Lemma~\ref{lemgddom}
shows that, for any fixed number of deletions over a finite horizon
$N$, the greedy deletion policy (Definition~\ref{defgdRule}) yields
the smallest area under the resulting sample path, $Q$, over $ \{
1,\ldots, N \}$. The main idea of proof is to show that the area
under $Q$ after applying $\pi_{\mathrm{NOB}}$ is asymptotically the same as
that of the greedy policy, as $N\to\infty$ and $\lambda\to1$ (in
this particular order of limits). In some sense, this means that the
jobs in $M^\Psi$ account for almost all of the delays in the system,
as $\lambda\to1$. The following technical lemma is useful.

%
\begin{lem}
\label{lemtopk}
For a finite set $S\subset\mathbb{R}$ and $k\in\mathbb{N}$, define
\[
f(S,k) = \frac{\mathrm{sum\ of\ the}\ k\ \mathrm{largest\ elements\ in}\ S}{|S|}.
\]
Let $ \{X_i\dvtx 1\leq i \leq n \}$ be i.i.d. random variables
taking values in $\mathbb{Z}_{+}$, where \mbox{$\mathbb{E} (X_1
)<\infty$}. Then for any sequence of random variables $ \{H_n\dvtx n\in
\mathbb{N} \}$, with\vadjust{\goodbreak} \mbox{$H_n \lesssim\alpha n$} a.s. as $n \to
\infty$ for some $\alpha\in(0,1)$, we have
%
%
\begin{equation}
\qquad \limsup_{n \to\infty} f \bigl( \{X_i\dvtx 1\leq i \leq n \},H_n \bigr) \leq\mathbb{E} \bigl(X_1\cdot\mathbb{I}
\bigl(X_1\geq\widebar{F}{}^{ -1}_{X_1}(\alpha)
\bigr) \bigr)\qquad\mbox{a.s.},
\end{equation}
where $\widebar{F}{}^{ -1}_{X_1}(y) = \min\{x\in\mathbb{N}\dvtx
\mathbb{P} (X_1\geq x )< y \}$.
\end{lem}

\begin{pf}
See Appendix~\ref{applemtopk}.
\end{pf}

Fix an initial sample path $Q^0$. We will denote by $M^\Psi= \{
m^\Psi_i\dvtx i\in\mathbb{N} \}$ the deletion sequence generated
by $\pi_{\mathrm{NOB}}$ on $Q^0$. Define
%
%
\begin{equation}
l (n ) = n-\max_{1\leq i \leq I (M^\Psi,n )} |E_i|, \label{eqldef}
\end{equation}
where $E_i$ is the $i$th deletion epoch of $M^\Psi$, defined in
equation~(\ref{eqEjdef}). Since $Q^0[n]\geq Q^0[m_i]$ for all $i\in
\mathbb{N}$, it is easy to check that
\[
\Delta_P \bigl(D \bigl(Q^0, \bigl\{m^\Psi_j\dvtx
1\leq j \leq i-1 \bigr\} \bigr),n,m^\Psi_i \bigr) =
n-m^\Psi_i +1
\]
for all $i\in\mathbb{N}$. The function $l$ was defined so that the
first $I(M^\Psi,l(n))$ deletions made by a greedy rule over the
horizon $ \{1,\ldots,n \}$ are exactly $ \{1,\ldots,
l(n) \}\cap M^\Psi$. More formally, we have the following lemma.

%
\begin{lem}
\label{lemGdoverlap}
Fix $n\in\mathbb{N}$, and let $M^G=G (Q^{0},n,I (M^\Psi,l (n ) ) )$. Then $m_{i}^{G}=m^\Psi_i$, for
all $i\in\{1,\ldots, I (M^\Psi,l(n) ) \}$.
\end{lem}

Fix $K\in\mathbb{N}$, and an arbitrary feasible deletion sequence,
$\widetilde{M}$, generated by a policy in $\Pi_\infty$. We can write
%
%
\begin{eqnarray}\label{eqIineq}
I \bigl(\widetilde{M}, m^\Psi_K \bigr)
&=& I \bigl(M^\Psi, l \bigl(m^\Psi_K \bigr)
\bigr) + \bigl(I \bigl(M^\Psi, m^\Psi_K \bigr) -
I \bigl(M^\Psi,l \bigl(m^\Psi_K \bigr) \bigr)
\bigr)\nonumber
\\
&&{} + \bigl(I \bigl(\widetilde{M}, m^\Psi_K \bigr)-I
\bigl(M^\Psi, m^\Psi_K \bigr) \bigr)\nonumber
\\
&=& I \bigl(M^\Psi, l \bigl(m^\Psi_K \bigr)
\bigr) + \bigl(K - I \bigl(M^\Psi,l \bigl(m^\Psi_K
\bigr) \bigr) \bigr)
\\
&&{}  + \bigl(I \bigl(\widetilde{M}, m^\Psi_K \bigr)-I
\bigl(M^\Psi, m^\Psi_K \bigr) \bigr)\nonumber
\\
&=& I \bigl(M^\Psi, l \bigl(m^\Psi_K \bigr)
\bigr) + h(K), \nonumber
\end{eqnarray}
where
%
%
\begin{equation}
h(K) = \bigl(K - I \bigl(M^\Psi,l \bigl(m^\Psi_K
\bigr) \bigr) \bigr) + \bigl(I \bigl(\widetilde{M}, m^\Psi_K
\bigr)-I \bigl(M^\Psi, m^\Psi_K \bigr) \bigr).
\end{equation}
We have the following characterization of $h$.

%
\begin{lem}
\label{lemhK}
$ h(K) \lesssim\frac{1-\lambda}{\lambda-(1-p)} \cdot K$, as $K\to
\infty$, a.s.
\end{lem}

\begin{pf}
See Appendix~\ref{applemhK}.\vadjust{\goodbreak}
\end{pf}

Let
%
%
\begin{equation}
M^{G,n} = G \bigl(Q^0,n, I (\widetilde{M},n ) \bigr),
\end{equation}
where the greedy deletion map $G$ was defined in Definition~\ref{defgdRule}. By Lemma~\ref{lemGdoverlap} and the definition of
$M^{G,n}$, we have that
%
%
\begin{equation}
M^\Psi\cap\bigl\{1,\ldots, l \bigl(m^\Psi_K
\bigr) \bigr\} \subset M^{G,m^\Psi_K}.
\end{equation}
Therefore, we can write
%
%
\begin{equation}
M^{G,m^\Psi_K} = \bigl(M^\Psi\cap\bigl\{1,\ldots, l
\bigl(m^\Psi_K \bigr) \bigr\} \bigr) \cup
\widebar{M}{}^G_K, \label{eqMGoverline}
\end{equation}
where $\widebar{M}{}^G_K\stackrel{\triangle}{=}M^{G,m^\Psi_K}
\setminus(M^\Psi\cap\{1,\ldots, l (m^\Psi
_K ) \} )$. Since $\llvert M^{G,m^\Psi_K}\rrvert =
I (\widetilde{M},m^\Psi_K )$ by definition, by
equation~(\ref
{eqIineq}),
%
\begin{equation}
\bigl\llvert\widebar{M}{}^G_K\bigr\rrvert= h(K).
\end{equation}

We have
%
%
\begin{eqnarray} \label{eqSdiff}
&& S \bigl(D \bigl(Q^0,M^\Psi\bigr),m^\Psi_K
\bigr) - S \bigl(D \bigl(Q^0,\widetilde{M} \bigr),m^\Psi_K
\bigr)\nonumber
\\
&&\qquad \stackrel{\mathrm{(a)}} {\leq}  S \bigl(D \bigl(Q^0,M^\Psi
\bigr),m^\Psi_K \bigr) - S \bigl(D \bigl(Q^0,M^{G,m^\Psi_K}
\bigr),m^\Psi_K \bigr)
\\
&&\qquad \stackrel{\mathrm{(b)}} {=} \Delta\bigl(D \bigl(Q^0,M^\Psi
\bigr),m^\Psi_K,\widebar{M}{}^G_K
\bigr),\nonumber
\end{eqnarray}
where (a) is based on the dominance of the greedy policy over any
finite horizon (Lemma~\ref{lemgddom}), and (b) follows from
equation~(\ref{eqMGoverline}).

Finally, we claim that there exists $g(x)\dvtx \mathbb{R}\to\mathbb
{R}_+$, with $g(x)\to0$ as $x\to1$, such that
%
%
\begin{equation}
\limsup_{K \to\infty} \frac{\Delta(D (Q^0,M^\Psi
),m^\Psi_K,\widebar{M}{}^G_K )}{m^\Psi_K} \leq g(\lambda)\qquad\mbox{a.s.}\label{eqdeltadim}
\end{equation}
Equations~(\ref{eqSdiff}) and (\ref{eqdeltadim}) combined imply that
%
%
\begin{eqnarray}
C (p,\lambda, \pi_{\mathrm{NOB}} ) &=& \limsup_{K\to\infty}
\frac
{S (D (Q^0,M^\Psi),m^\Psi_K )}{m^\Psi_K}
\nonumber
\\
&\leq& g(\lambda)+ \limsup_{K\to\infty} \frac{S (D
(Q^0,\widetilde{M} ),m^\Psi_K )}{m^\Psi_K},
\\
&= & g(\lambda)+ \limsup_{n \to\infty} \frac{S (D
(Q^0,\widetilde{M} ),n )}{n}\qquad\mbox{a.s.},\nonumber
\end{eqnarray}
which shows that
\[
C (p,\lambda, \pi_{\mathrm{NOB}} ) \leq g(\lambda) + \inf_{\pi
\in\Pi_\infty}C
(p,\lambda, \pi).
\]
Since $g(\lambda)\to0$ as $\lambda\to1$, this proves Proposition
\ref{propnobopt}.

To show equation~(\ref{eqdeltadim}), denote by $Q$ the sample path
after applying $\pi_{\mathrm{NOB}}$,
\[
Q = D \bigl(Q^0,M^\Psi\bigr)\vadjust{\goodbreak}
\]
and by $V_i$ the area under $Q$ within $E_i$,
\[
V_i = \sum_{n=m^\Psi_i}^{m^\Psi_{i+1}-1} Q [n ].
\]
An example of $V_i$ is illustrated as the area of the shaded region in
Figure~\ref{figwater2}. By Proposition~\ref{propQRW}, $Q$ is a
Markov chain, and so is the process $W[n]= (Q[n], Q[n+1] )$.
By Lemma~\ref{lemqmi0}, $E_i$ corresponds to the indices between two
adjacent returns of the chain $W$ to state $(0,0)$. Since the $i$th
return of a Markov chain to a particular state is a stopping time, it
can be shown, using the strong Markov property of $W$, that the
segments of $Q$, $ \{Q[n]\dvtx n\in E_i \}$, are mutually
independent and identically distributed among different values of $i$.
Therefore, the $V_i$'s are i.i.d. Furthermore,
%
%
\begin{equation}
\mathbb{E} (V_1 )\stackrel{\mathrm{(a)}} {\leq} \mathbb{E}
\bigl(|E_1|^2 \bigr) \stackrel{\mathrm{(b)}} {<} \infty,
\end{equation}
where (a) follows from the fact that $|Q[n+1]-Q[n]|\leq1$ for all
$n$, and hence $V_i\leq|E_i|^2$ for any sample path of $Q^0$, and
(b) from the exponential tail bound on $\mathbb{P}(|E_1|\geq x)$,
given in equation~(\ref{eqEiExp}).

Since the value of $Q$ on the two ends of $E_i$, $m^\Psi_i$ and
$m^\Psi_{i+1}-1$, are both zero, each additional deletion within
$E_i$ cannot produce a marginal decrease of area under $Q$ of more than
$V_i$; cf. Figure~\ref{figwater2}. Therefore, the\vspace*{1pt} value of~$\Delta
(D (Q^0,M^\Psi),m^\Psi_K,\widebar{M}{}^G_K )$
can be no greater than the sum of the $h(K)$ largest $V_i$'s over the
horizon $n\in\{1,\ldots, m^\Psi_K \}$. We have
%
%
\begin{eqnarray}
&& \limsup_{K \to\infty} \frac{\Delta(D (Q^0,M^\Psi
),m^\Psi_K,\widebar{M}{}^G_K )}{m^\Psi_K}
\nonumber
\\
&&\qquad = \limsup_{K \to\infty} f \bigl( \{V_i\dvtx  1\leq i \leq K
\},h(K) \bigr) \cdot\frac{K}{m^\Psi_K}
\nonumber\\[-8pt]\\[-8pt]
&&\qquad \stackrel{\mathrm{(a)}} {=}  \limsup_{K \to\infty} f \bigl(
\{V_i\dvtx  1\leq i \leq K \},h(K) \bigr) \cdot\frac{\lambda+1-p}{\lambda
-(1-q)}
\nonumber
\\
&&\qquad \stackrel{\mathrm{(b)}} {=}  \mathbb{E} \biggl(V_1\cdot\mathbb{I}
\biggl(X_1\geq\widebar{F}{}^{ -1}_{V_1} \biggl(
\frac{1-\lambda}{\lambda
-(1-p)} \biggr) \biggr) \biggr) \cdot\frac{\lambda+1-p}{\lambda-(1-q)},\nonumber
\end{eqnarray}
where (a) follows from equation~(\ref{eqmilim}), and (b) from
Lemmas~\ref{lemtopk} and~\ref{lemhK}. Since $\mathbb{E}
(V_1 )< \infty$, and $\widebar{F}{}^{ -1}_{V_1}(x)\to\infty$
as $x \to0$, it follows that
\[
\mathbb{E} \biggl(V_1\cdot\mathbb{I} \biggl(X_1\geq
\widebar{F}{}^{
-1}_{V_1} \biggl(\frac{1-\lambda}{\lambda-(1-p)} \biggr)
\biggr) \biggr) \to0
\]
as $\lambda\to1$. Equation~(\ref{eqdeltadim}) is proved by setting
\[
g(\lambda) = \mathbb{E} \biggl(V_1\cdot\mathbb{I} \biggl(X_1\geq
\widebar{F}{}^{ -1}_{V_1} \biggl(\frac{1-\lambda}{\lambda-(1-p)}
\biggr) \biggr) \biggr) \cdot\frac{\lambda+1-p}{\lambda-(1-q)}.
\]
This completes the proof of Proposition~\ref{propnobopt}.
\end{pf*}

\subsubsection{\texorpdfstring{Why not use greedy?}{Why not use greedy}}\label{sec6.3.1} The proof of Proposition~\ref{propnobopt} relies on a sample-path-wise coupling to the performance
of a greedy deletion rule. It is then only natural to ask: since the
time horizon is indeed finite in all practical applications, why do not
we simply use the greedy rule as the preferred offline policy, as
opposed to $\pi_{\mathrm{NOB}}$?

There are at least two reasons for focusing on $\pi_{\mathrm{NOB}}$ instead of
the greedy rule. First, the structure of the greedy rule is highly
global, in the sense that each deletion decision uses information of
the entire sample path over the horizon. As a result, the greedy rule
tells us little on how to design a good policy with a \emph{fixed}
lookahead window (e.g., Theorem~\ref{teolookahead}). In contrast, the
performance analysis of $\pi_{\mathrm{NOB}}$ in Section~\ref
{secperformoffline} reveals a highly \emph{regenerative} structure:
the deletions made by $\pi_{\mathrm{NOB}}$ essentially depend only on the
dynamics of $Q^0$ in the same deletion epoch (the $E_i$'s), and what
happens beyond the current epoch becomes irrelevant. This is the key
intuition that led to our construction of the finite-lookahead policy
in Theorem~\ref{teolookahead}. A second (and perhaps minor) reason is
that of computational complexity. By a small sacrifice in performance,
$\pi_{\mathrm{NOB}}$ can be efficiently implemented using a linear-time
algorithm (Section~\ref{seclinearalgo}), while it is easy to see
that a naive implementation of the greedy rule would require
super-linear complexity with respect to the length of the horizon.

\subsection{\texorpdfstring{Proof of Theorem \protect\ref{teooffline}.}{Proof of Theorem 2}}\label{secpfthmoffline}

The fact that $\pi_{\mathrm{NOB}}$ is feasible follows from equation~(\ref
{eqIMnlim}) in Lemma~\ref{lemnobbasic}, that is,
\[
\limsup_{n\to\infty}\frac{1}{n}I \bigl(M^\Psi,n
\bigr) \leq\frac{\lambda- (1-p )}{\lambda+1-p}<\frac{p}{\lambda
+1-p}\qquad\mbox{a.s.}
\]

Let $ \{\widetilde{Q}[n]\dvtx n\in\mathbb{Z}_{+} \}$ be the resulting
sample path after applying $\pi_{\mathrm{NOB}}$ to the initial sample path
$ \{Q^0[n]\dvtx n\in\mathbb{Z}_{+} \}$, and let
\[
Q[n]=\widetilde{Q} \bigl[n+m^\Psi_1 \bigr]\qquad\forall
n \in\mathbb{N},
\]
where $m^\Psi_1$ is the index of the first deletion made by $\pi_{\mathrm{NOB}}$.
Since $\lambda>1-p$, the random walk $Q^0$ is transient, and hence
$m^\Psi_1<\infty$ almost surely. We have that, almost surely,
%
%
\begin{eqnarray}\label{eqcplam}
C (p,\lambda, \pi_{\mathrm{NOB}} ) &=& \lim_{n\to\infty}
\frac
{1}{n}\sum_{i=1}^n
\widetilde{Q}[i]
\nonumber
\\
&=& \lim_{n\to\infty}\frac{1}{n}\sum
_{i=1}^{m^\Psi_1} \widetilde{Q}[i] + \lim
_{n\to\infty}\frac{1}{n}\sum_{i=1}^n
Q[i]
\\
&=& \frac{1-p}{\lambda-(1-p)}, \nonumber
\end{eqnarray}
where the last equality follows from equation~(\ref{eqQavg}) in
Proposition~\ref{proppnobperf}, and the fact that $m_1<\infty$
almost surely. Letting $\lambda\to1$ in equation~(\ref{eqcplam})
yields the finite limit of delay under heavy traffic,
\[
\lim_{\lambda\to1}C (p,\lambda, \pi_{\mathrm{NOB}} ) = \lim
_{\lambda\to1} \frac{1-p}{\lambda-(1-p)} = \frac{1-p}{p}.
\]

Finally, the delay optimality of $\pi_{\mathrm{NOB}}$ in heavy traffic was
proved in Proposition~\ref{propnobopt}, that is, that
\[
\lim_{\lambda\to1}C (p,\lambda,\pi_{\mathrm{NOB}} )= \lim
_{\lambda\to1}C^*_{\Pi_\infty} (p,\lambda).
\]
This completes the proof of Theorem~\ref{teooffline}.

\section{\texorpdfstring{Policies with a finite lookahead.}{Policies with a finite lookahead}}\label{sec7}
\label{secfinitelookahead}

\subsection{\texorpdfstring{Proof of Theorem \protect\ref{teolookahead}.}{Proof of Theorem 3}}\label{secpfthmlookahead}

As pointed out in the
discussion preceding Theorem~\ref{teolookahead}, for any initial
sample path and $w<\infty$, an arrival that is deleted under the\vadjust{\goodbreak} $\pi
_{\mathrm{NOB}}$ policy will also be deleted under $\pi_{\mathrm{NOB}}^{w}$. Therefore,
the delay guarantee for $\pi_{\mathrm{NOB}}$ (Theorem~\ref{teooffline})
carries over to $\pi_{\mathrm{NOB}}^{w(\lambda)}$, and for the rest of the
proof, we will be focusing on showing that $\pi_{\mathrm{NOB}}^{w(\lambda)}$
is feasible under an appropriate scaling of $w(\lambda)$. We begin by
stating an exponential tail bound on the distribution of the
discrete-time predictive window, $W(\lambda,n)$, defined in
equation~(\ref{eqwdiscrete}),
\[
W(\lambda,n) = \max\bigl\{ k \in\mathbb{Z}_{+}\dvtx  T_{n+k}
\leq T_{n}+w(\lambda) \bigr\}.
\]
It is easy to see that $ \{W (\lambda,m^\Psi_i )\dvtx i
\in\mathbb{N} \}$ are i.i.d., with $W (\lambda,m^\Psi
_1 )$ distributed as a Poisson random variable with mean
$(\lambda+1-p)w(\lambda)$. Since
\[
\mathbb{P} \bigl(W \bigl(\lambda,m^\Psi_1 \bigr) \geq x
\bigr) \leq\mathbb{P} \Biggl( \sum_{k=1}^{ \lfloor w(\lambda)
\rfloor}
X_k \Biggr),
\]
where the $X_k$ are i.i.d. Poisson random variables with mean $\lambda
+(1-p)$, applying the Chernoff bound, we have that, there exist $c,d>0$
such that
%
%
\begin{equation}
\mathbb{P} \biggl(W \bigl(\lambda,m^\Psi_1 \bigr) \geq
\frac
{\lambda+1-p}{2} \cdot w(\lambda) \biggr) \leq c\cdot\exp\bigl(-d\cdot w(
\lambda)\bigr)
\end{equation}
for all $w(\lambda)>0$.

We now analyze the deletion rate resulted by the $\pi_{\mathrm{NOB}}^{w(\lambda
)}$ policy. For the pure purpose of analysis (as opposed to practical
efficiency), we will consider a new deletion policy, denoted by $\sigma
^{w(\lambda)}$, which can be viewed as a relaxation of $\pi
_{\mathrm{NOB}}^{w(\lambda)}$.

%
\begin{defn}
Fix $w\in\mathbb{R}_{+}$. The deletion policy $\sigma^w$ is defined
such that for each deletion epoch $E_i$, $i\in\mathbb{N}$:
\begin{longlist}[(2)]
\item[(1)] if $|E_i| \leq W (\lambda,m^\Psi_i )$, then
only the first arrival of this epoch, namely, the arrival in slot
$m^\Psi_i$, is deleted;
\item[(2)] otherwise, all arrivals within this epoch are deleted.
\end{longlist}
\end{defn}

It is easy to verify that $\sigma^{w}$ can be implemented with $w$
units of look-ahead, and the set of deletions made by $\sigma
^{w(\lambda)}$ is a strict superset of $\pi_{\mathrm{NOB}}^{w(\lambda)}$
almost surely. Hence, the feasibility of $\sigma^{w(\lambda)}$ will
imply that of $\pi_{\mathrm{NOB}}^{w(\lambda)}$.

Denote by $D_i$ the number of deletions made by $\sigma^{w(\lambda)}$
in the epoch $E_i$. By the construction of the policy, the $D_i$ are
i.i.d., and depend only on the length of $E_i$ and the number of
arrivals within. We have\footnote{For simplicity of notation, we
assume that $\frac{\lambda+1-p}{2}\cdot w(\lambda)$ is always an
integer. This does not change the scaling behavior of $w(\lambda)$. }
%
%
\begin{eqnarray}\label{eqED1}
\mathbb{E} (D_1 )&\leq& 1+ \mathbb{E} \bigl[|E_i| \cdot\mathbb{I}
\bigl(|E_i|\geq W \bigl(\lambda, m^\Psi_i \bigr)
\bigr) \bigr]
\nonumber
\\
&\leq& 1+ \mathbb{E} \biggl[|E_i|\cdot\mathbb{I}
\biggl(|E_i| \geq\frac{\lambda+1-p}{2}\cdot w(\lambda) \biggr) \biggr]
\nonumber
\\
&&{} + \mathbb{E} \bigl(|E_i| \bigr)\cdot\mathbb{P} \biggl( W \bigl(
\lambda,m^\Psi_i \bigr) \leq\frac{\lambda+1-p}{2}\cdot w(
\lambda) \biggr)
\nonumber\\[-8pt]\\[-8pt]
&\leq& 1 + \Biggl(\sum_{k = ((\lambda+1-p)/2)\cdot w(\lambda
)}^{\infty} k \cdot a
\cdot\exp(-b\cdot k) \Biggr)\nonumber
\\
&&{} + \frac{\lambda}{\lambda-(1-p)}\cdot c\cdot\exp\bigl(-d\cdot w(
\lambda)\bigr)\nonumber
\\
&\stackrel{\mathrm{(a)}} {\leq}&  1+ h \cdot w(\lambda) \cdot\exp\bigl(-l\cdot
w(\lambda)\bigr)\nonumber
\end{eqnarray}
for some $h,l>0$, where (a) follows from the fact that $\sum
_{k=n}^\infty k\cdot\exp(-b\cdot k) = \mathcal{O} (n \cdot\exp
(-b \cdot n) )$ as $n\to\infty$.

Since the $D_i$ are i.i.d., using basic renewal theory, it is not
difficult to show that the average rate of deletion in discrete time
under the policy $\sigma^{w(\lambda)}$ is equal to $\frac{\mathbb
{E}(D_1)}{\mathbb{E}(E_1)}$. In order for the policy to be feasible,
one must have that
%
%
\begin{equation}
\frac{\mathbb{E}(D_1)}{\mathbb{E}(E_1)} = \frac{\mathbb
{E}(D_1)}{\lambda} \leq\frac{p}{\lambda+ 1-p}. \label{eqdelrate1}
\end{equation}
By equations~(\ref{eqED1}) and (\ref{eqdelrate1}), we want to
ensure that
\[
\frac{p\lambda}{\lambda- (1-p)} \geq1+h\cdot w(\lambda)\cdot\exp\bigl
(-l\cdot w(\lambda)
\bigr),
\]
which yields, after taking the logarithm on both sides,
%
%
\begin{equation}
w(\lambda)\geq\frac{1}{b} \log\biggl(\frac{1}{1-\lambda} \biggr) +
\frac{1}{b} \log\biggl(\frac{[\lambda-(1-p)]\cdot h\cdot
w(\lambda)}{1-p} \biggr).
\end{equation}
It is not difficult to verify that for all $p\in(0,1)$ there exists a
constant $C$ such that the above inequality holds for all $\lambda\in
(1-p,1)$, by letting $w(\lambda)=C\log(\frac{1}{1-\lambda})$. This
proves the feasibility of $\sigma^{w(\lambda)}$, which implies that
$\pi_{\mathrm{NOB}}^{w(\lambda)}$ is also feasible. This completes the proof
of Theorem~\ref{teolookahead}.

\section{\texorpdfstring{Concluding remarks and future work.}{Concluding remarks and future work}}\label{sec8}
\label{secconclusions}

The main objective of this paper is to study the impact of future
information on the performance of a class of admissions control
problems, with a constraint on the time-average rate of redirection.
Our model is motivated as a study of a dynamic resource allocation
problem between slow (congestion-prone) and fast (congestion-free)
processing resources. It could also serve as a simple canonical model
for analyzing delays in large server farms or cloud clusters with
resource pooling \cite{TX12}; cf. Appendix~\ref{secresourcepooling}. Our main results show that the availability of
future information can dramatically reduce the delay experienced by
admitted customer: the delay converges to a finite constant even as the
traffic load approaches the system capacity (``heavy-traffic delay
collapse''), if the decision maker is allowed for a sufficiently large
lookahead window (Theorem~\ref{teolookahead}).

There are several interesting directions for future exploration. On the
theoretical end, a main open question is whether a matching lower-bound
on the amount of future information required to achieve the
heavy-traffic delay collapse can be proved (Conjecture~\ref{conjinfolowerbound}), which, together with the upper bound given in
Theorem~\ref{teolookahead}, would imply a duality between delay and
the length of lookahead into the future.

Second, we believe that our results can be generalized to the cases
where the arrival and service processes are non-Poisson. We note that
the $\pi_{\mathrm{NOB}}$ policy is indeed feasible for a wide range of
non-Poisson arrival and service processes (e.g., renewal processes), as
long as they satisfy a form of strong law of large number, with
appropriate time-average rates (Lemma~\ref{lemnobbasic}). It seems
more challenging to generalize results on the optimality of $\pi
_{\mathrm{NOB}}$ and the performance guarantees. However, it may be possible to
establish a generalization of the delay optimality result using
limiting theorems (e.g., diffusion approximations). For instance, with
sufficiently well-behaved arrival and service processes, we expect that
one can establish a result similar to Proposition~\ref{propQRW} by
characterizing the resulting queue length process from $\pi_{\mathrm{NOB}}$ as
a reflected Brownian motion in $\mathbb{R}_+$, in the limit of
$\lambda\to1$ and $p\to0$, with appropriate scaling.

Another interesting variation of our problem is the setting where each
job comes with a prescribed \emph{size}, or \emph{workload}, and the
decision maker is able to observe both the arrival times and workloads
of jobs up to a finite lookahead window. It is conceivable that many
analogous results can be established for this setting, by studying the
associated workload (as opposed to queue length) process, while the
analysis may be less clean due to the lack of a simple
random-walk-based description of the system dynamics. Moreover, the
\emph{server} could potentially exploit additional information of the
jobs' workloads in making scheduling decisions, and it is unclear what
the performance and fairness implications are for the design of
admissions control policies.

There are other issues that need to be addressed if our offline
policies (or policies with a finite lookahead) are to be applied in
practice. A most important question can be the impact of \emph
{observational noise} to performance, since in reality the future seen
in the lookahead window cannot be expected to match the actual
realization exactly. We conjecture, based on the analysis of $\pi
_{\mathrm{NOB}}$, that the performance of both $\pi_{\mathrm{NOB}}$, and its
finite-lookahead version, is robust to small noises or perturbations
(e.g., if the actual sample path is at most $\varepsilon$ away from the
predicted one), while it remains to thoroughly verify and quantify the
extend of the impact, either empirically or through theory. Also, it is
unclear what the best practices should be when the lookahead window is
very small relative to the traffic intensity $\lambda$ ($w \ll\log
{\frac{1}{1-\lambda}}$), and this\vspace*{1pt} regime is not covered by the
results in this paper (as illustrated in Figure~\ref{figcollapse}).

\begin{appendix}
\section{\texorpdfstring{Additional Proofs}{Additional Proofs}}\label{secA}

\subsection{\texorpdfstring{Proof of Lemma \protect\ref{lemnobbasic}.}{Proof of Lemma 2}}\label{applemnobbasic}
Since $\lambda> 1-p$, with probability one, there exists $T<\infty$
such that the continuous-time queue length process without deletion
satisfies $Q^0(t)>0$ for all $t\geq T$. Therefore, without any
deletion, all service tokens are matched with some job after time $T$.
By the stack interpretation, $\pi_{\mathrm{NOB}}$ only deletes jobs that would
not have been served, and hence does not change the original matching
of service tokens to jobs. This proves the first claim.

By the first claim, since all subsequent service tokens are matched
with a job after some time $T$, there exists some $N<\infty$, such that
%
%
\begin{equation}
\widetilde{Q}[n] = \widetilde{Q}[N] + \bigl(A[n]-A[N] \bigr)- \bigl(S[n]-S[N]
\bigr) - I \bigl(M^\Psi,n \bigr) \label{eqqevolv}
\end{equation}
for all $n \geq N$, where $A[n]$ and $S[n]$ are the cumulative numbers
of arrival and service tokens by slot $n$, respectively. The\vspace*{1pt} second
claim follows by multiplying both sides of equation~(\ref{eqqevolv})
by $\frac{1}{n}$, and using the fact that $\lim_{n\to\infty}\frac
{1}{n}A[n]=\frac{\lambda}{\lambda+1-p}$ and $\lim_{n\to\infty
}\frac{1}{n}S[n]=\frac{1-p}{\lambda+1-p}$ a.s., $\widetilde
{Q}[n]\geq0$
for all $n$ and $\widetilde{Q}[N]<\infty$ a.s.

\subsection{\texorpdfstring{Proof of Lemma \protect\ref{lemqmi0}.}{Proof of Lemma 4}}\label{secA.2}
\label{applemqmi0}

(1)~Recall\vspace*{1pt} the point-wise deletion map, $D_P (Q,n
)$, defined in Definition~\ref{defdelMap}. For any initial sample
path $Q^0$, let $Q^1 = D_P(Q^0,m)$ for some $m \in\mathbb{N}$. It is
easy to see that, for all $n>m$, $Q^1[n] = Q^0[n]-1$, if and only if
$Q^0[s]\geq1$ for all $s \in\{m+1,\ldots, n \}$.
Repeating this argument $I(M,n)$ times, we have that
%
%
\begin{equation}
Q[n] = \widetilde{Q}[n+m_1] = Q^0[n+m_1] -
I (M,n+m_1 ), \label{eqQnevol1}
\end{equation}
if and only if for all $k \in\{1,\ldots,I(M,n+m_1) \}$,
%
%
\begin{equation}
Q^0[s] \geq k\qquad\mbox{for all $s\in\{m_k+1,\ldots,n+m_1 \}$}. \label{eqqmi0}
\end{equation}
Note that equation~(\ref{eqqmi0}) is implied by (and in fact,
equivalent to) the definition of the $m_k$'s (Definition~\ref{defnob}), namely, that for all $k\in\mathbb{N}$, $Q^0[s]\geq k$ for
all $s\geq m_k+1$. This proves the first claim.

(2)~Suppose $Q[n]=Q[n-1]=0$. Since $\mathbb{P} (Q^0[t]
\neq Q^0[t-1] \mid Q^0[t-1]>0 )=1$ for all $t\in\mathbb{N}$
[cf. equation~(\ref{eqQ0trans1})] at least one deletion occurs on the
slots $ \{n-1+m_1,n+m_1 \}$. If the deletion occurs on
$n+m_1$, we are done. Suppose a deletion occurs on $n-1+m_1$.
Then\vadjust{\goodbreak}
$Q^0[n+m_1]\geq Q^0[n-1+m_1]$, and hence
\[
Q^0[n+m_1]= Q^0[n-1+m_1]+1,
\]
which implies that a deletion must also occur on $n+m_1$, for otherwise
$Q[n]=Q[n-1]+1 = 1 \neq0$. This shows that $n=m_i-m_1$ for some $i\in
\mathbb{N}$.

Now, suppose that $n=m_i-m_1$ for some $i\in\mathbb{N}$. Let
%
%
\begin{equation}
n_k = \inf\bigl\{n\in\mathbb{N}\dvtx  Q^0[n]=k\mbox{ and }
Q^0[t]\geq k,\ \forall t\geq n \bigr\}.
\end{equation}
Since the random walk $Q^0$ is transient, and the magnitude of its step
size is at most~$1$, it follows that $n_k<\infty$ for all $k\in
\mathbb{N}$ a.s., and that $m_k= n_k$, $\forall k\in\mathbb{N}$. We have
%
%
\begin{eqnarray} \label{eqQnevol2}
Q[n]
&\stackrel{\mathrm{(a)}} {=} & Q^0 [n+m_1]-I (M,n+m_1)
\nonumber
\\
&= & Q^0[m_i] - I (M,m_i )
\nonumber\\[-8pt]\\[-8pt]
&\stackrel{\mathrm{(b)}} {=} & Q^0[n_i] - i\nonumber
\\
&= & 0,\nonumber
\end{eqnarray}
where (a) follows from equation~(\ref{eqQnevol1}) and (b) from
the fact that $n_i=m_i$. To show that $Q[n-1]=0$, note that since
$n=m_i-m_1$, an arrival must have occurred in $Q^0$ on slot $m_i$, and
hence $Q^0[n-1+m_1]=Q^0[n+m_1]-1$. Therefore, by the definition of $m_i$,
\[
Q^0[t] -Q^0[n-1+m_1]= \bigl(Q^0[t]-Q^0[n+m_1]
\bigr)+1 \geq0\qquad\forall t\geq n+m_1,
\]
which implies that $n-1 = m_{i-1}-m_1$, and hence $Q[n-1]=0$, in light
of equation~(\ref{eqQnevol2}). This proves the claim.

(3)~For all $n\in\mathbb{Z}_{+}$, we have
%
%
\begin{eqnarray}
Q[n] &=& Q [m_{I (M,n+m_1 )}-m_1 ]+ \bigl(Q[n]-Q
[m_{I (M,n+m_1 )}-m_1 ] \bigr)
\nonumber
\\
&\stackrel{\mathrm{(a)}} {=}& Q[n]-Q [m_{I (M,n+m_1 )}-m_1 ]
\nonumber\\[-8pt]\\[-8pt]
&\stackrel{\mathrm{(b)}} {=}& Q^0[n+m_1]-Q^0
[m_{I (M,n+m_1
)} ]\nonumber
\\
&\stackrel{\mathrm{(c)}} {=} & 0,\nonumber
\end{eqnarray}
where (a) follows from the second claim [cf. equation~(\ref
{eqQnevol2})], (b) from the fact that there is no deletion on any
slot in $ \{I (M,n+m_1 ),\ldots, n+m_1 \}$ and
(c) from the fact that $n+m_1\geq I (M,n+m_1 )$ and
equation~(\ref{eqmiinc}).
%

\subsection{\texorpdfstring{Proof of Lemma \protect\ref{lemdualRW}.}{Proof of Lemma 5}}\label{applemdualRW}
Since the random walk $X$ lives in $\mathbb{Z}_{+}$ and can take jumps
of size at most $1$, it suffices to verify that
\[
\mathbb{P} \Bigl(X[n+1]=x_1+1 | X[n]=x_1, \min
_{r \geq n+1} X[r] = 0 \Bigr) = 1-q\vadjust{\goodbreak}
\]
for all $x_1\in\mathbb{Z}_{+}$. We have
%
%
\begin{eqnarray} \label{eqXcond1}
&& \mathbb{P} \Bigl(X[n+1]=x_1+1 | X[n]=x_1, \min
_{r \geq n+1} X[r] = 0 \Bigr)
\nonumber
\\
&&\qquad = \mathbb{P} \Bigl(X[n+1]=x_1+1, \min_{r \geq n+1} X[r] = 0
| X[n]=x_1 \Bigr)\nonumber
\\
&&\quad\qquad{} \bigm/
\mathbb{P} \Bigl(\min_{r \geq n+1} X[r] =
0 | X[n]=x_1 \Bigr)
\nonumber
\\
&&\qquad \stackrel{\mathrm{(a)}} {=} \Bigl(\mathbb{P} \bigl(X[n+1]=x_1+1 |
X[n]=x_1 \bigr)
\\
&&\hspace*{36pt}{}\times \mathbb{P} \Bigl(\min_{r \geq n+1} X[r] = 0
| X[n+1]=x_1+1 \Bigr)\Bigr)\nonumber
\\
&&\quad\qquad{}\bigm / \mathbb{P} \Bigl(\min_{r \geq n+1}X[r] = 0 | X[n]=x_1 \Bigr)\nonumber
\\
&&\qquad \stackrel{\mathrm{(b)}} {=} q\cdot\frac{h(x_1+1)}{h(x_1)},\nonumber
\end{eqnarray}
%
%
where
\[
h(x) = \mathbb{P} \Bigl(\min_{r \geq2} X[r] = 0 | X[1]=x \Bigr)
\]
and steps (a) and (b) follow from the Markov property and
stationarity of $X$, respectively. The values of $ \{h(x)\dvtx x\in
\mathbb{Z}_{+} \}$ satisfy the set of harmonic equations
%
%
\begin{equation}
\label{eqharm1} h(x)=\cases{ q\cdot h(x+1)+(1-q)\cdot h (x-1), &\quad
$x\geq1$,
\vspace*{5pt}\cr
q\cdot h(1)+1-q, &\quad$x=0$}
\end{equation}
with the boundary condition
%
%
\begin{equation}
\lim_{x\to\infty} h(x) =0. \label{eqharm2}
\end{equation}
Solving equations~(\ref{eqharm1}) and (\ref{eqharm2}), we obtain
the unique solution
\[
h(x) = \biggl(\frac{1-q}{q} \biggr)^x
\]
for all $x\in\mathbb{Z}_{+}$. By equation~(\ref{eqXcond1}), this
implies that
\[
\mathbb{P} \Bigl(X[n+1]=x_1+1 | X[n]=x_1, \min
_{r \geq n+1} X[r] = 0 \Bigr) = q\cdot\frac{1-q}{q}=1-q,
\]
which proves the claim.

\subsection{\texorpdfstring{Proof of Lemma \protect\ref{lemtopk}.}{Proof of Lemma 8}}\label{secA.4}
\label{applemtopk}
By the definition of $\widebar{F}{}^{ -1}_{X_1}$ and the strong law of
large numbers (SLLN), we have
%
%
\begin{equation}\quad
\lim_{n \to\infty} \frac{1}{n}\sum
_{i=1}^n \mathbb{I} \bigl(X_i \geq
\widebar{F}{}^{ -1}_{X_1}(\alpha) \bigr) = \mathbb{E} \bigl(
\mathbb{I} \bigl(X_i \geq\widebar{F}{}^{ -1}_{X_1}(
\alpha) \bigr) \bigr) < \alpha\qquad\mbox{a.s.} \label{eqslln1}\vadjust{\goodbreak}
\end{equation}
Denote by $S_{n,k}$ set of top $k$ elements in $ \{X_i\dvtx 1\leq i
\leq n \}$. By equation~(\ref{eqslln1}) and the fact that $H_n
\lesssim\alpha n$ a.s., there exists $N>0$ such that
\[
\mathbb{P} \bigl\{ \exists N \mbox{ s.t. } \min{S_{n,H_n}} \geq
\widebar{F}{}^{ -1}_{X_1}(\alpha),\ \forall n\geq N \bigr\} = 1,
\]
which implies that
%
%
\begin{eqnarray}
&& \limsup_{n \to\infty} f \bigl( \{X_i\dvtx 1\leq i \leq n \},H_n \bigr)
\nonumber
\\
&&\qquad \leq \limsup_{n\to\infty} \frac{1}{n}\sum
_{i=1}^n X_i\cdot\mathbb{I}
\bigl(X_i \geq\widebar{F}{}^{ -1}_{X_1}(\alpha)
\bigr)
\\
&&\qquad = \mathbb{E} \bigl(X_1\cdot\mathbb{I} \bigl(X_1\geq
\widebar{F}{}^{ -1}_{X_1}(\alpha) \bigr) \bigr)\qquad\mbox{a.s.},\nonumber
\end{eqnarray}
where the last equality follows from the SLLN. This proves our claim.

\subsection{\texorpdfstring{Proof of Lemma \protect\ref{lemhK}.}{Proof of Lemma 10}}\label{applemhK}
We begin by stating the following fact:
%
\begin{lem}
\label{lemmaxexpo}
Let $ \{X_i\dvtx i\in\mathbb{N} \}$ be i.i.d. random variables
taking values in $\mathbb{R}_+$, such that for some $a,b>0$, $\mathbb
{P} (X_1\geq x )\leq a\cdot\exp(-b\cdot x)$ for all $x\geq
0$. Then
\[
\max_{1\leq i \leq n} X_i = o(n)\qquad\mbox{a.s.}
\]
as $n\to\infty$.
\end{lem}
\begin{pf}
%
%
\begin{eqnarray}
\lim_{n\to\infty} \mathbb{P} \biggl(\max_{1\leq i \leq n}
X_i \leq\frac{2}{b}\ln n \biggr) &=& \lim_{n\to\infty}
\mathbb{P} \biggl(X_1 \leq\frac{2}{b}\ln n
\biggr)^n
\nonumber
\\
&\leq&\lim_{n\to\infty} \bigl(1-a\cdot\exp(- 2\ln n)
\bigr)^n
\nonumber\\[-8pt]\\[-8pt]
& =& \lim_{n\to\infty} \biggl(1-\frac{a}{n^2}
\biggr)^n
\nonumber
\\
& =& 1.\nonumber
\end{eqnarray}
In other words, $\max_{1\leq i \leq n} X_i \leq\frac{2}{b}\ln n$
a.s. as $n\to\infty$, which proves the claim.
\end{pf}

Since the $|E_i|$'s are i.i.d. with $\mathbb{E} (|E_1|
)=\frac{\lambda+1-p}{\lambda-(1-p)}$ (Proposition~\ref{proppnobperf}), we have that, almost surely,
%
%
\begin{equation}
\qquad m^\Psi_K = \sum_{i=0}^{K-1}
|E_i|\sim\mathbb{E} \bigl(|E_1| \bigr)\cdot K = \frac{\lambda+1-p}{\lambda-(1-p)}
\cdot K\qquad\mbox{as $K \to\infty$} \label{eqmnbKasymp1}
\end{equation}
by the strong law of large numbers.
By Lemma~\ref{lemmaxexpo} and equations~(\ref{eqEiExp}), we have
%
%
\begin{equation}
\max_{1\leq i \leq K}|E_i| = o(K)\qquad\mbox{a.s.}
\label{eqEsmallK}\vadjust{\goodbreak}
\end{equation}
as $K\to\infty$. By equation~(\ref{eqEsmallK}) and the fact that
$I (M^\Psi, m^\Psi_K )=K$, we have
%
%
\begin{eqnarray} \label{eqsmalloK1}
K - I \bigl(M^\Psi,l \bigl(m^\Psi_K \bigr)
\bigr) &=& K - I \Bigl(M^\Psi,m^\Psi_K - \max
_{1\leq i \leq K}|E_i| \Bigr)
\nonumber
\\
&\stackrel{\mathrm{(a)}} {\leq}& K - I \bigl(M^\Psi,m^\Psi_K
\bigr) + \max_{1\leq i \leq K}|E_i|
\nonumber\\[-8pt]\\[-8pt]
&=& \max_{1\leq i \leq K}|E_i|
\nonumber
\\
&=& o (K )\qquad\mbox{a.s.}\nonumber
\end{eqnarray}
as $K\to\infty$, where (a) follows from the fact that at most one
deletion can occur in a single slot, and hence $I(M,n+m)\leq I(M,n)+m$
for all $m,n\in\mathbb{N}$. Since $\widetilde{M}$ is feasible,
%
%
\begin{equation}
I (\widetilde{M},n ) \lesssim\frac{p}{\lambda+1-p}\cdot n \label{eqtilMratio}
\end{equation}
as $n\to\infty$. We have
\begin{eqnarray}
h(K)&=& \bigl(K - I \bigl(M^\Psi,l \bigl(m^\Psi_K
\bigr) \bigr) \bigr) + \bigl(I \bigl(\widetilde{M}, m^\Psi_K
\bigr)-I \bigl(M^\Psi, m^\Psi_K \bigr) \bigr)
\nonumber
\\
&\stackrel{\mathrm{(a)}} {\lesssim}& \bigl(K - I \bigl(M^\Psi,l
\bigl(m^\Psi_K \bigr) \bigr) \bigr) + \frac{p}{\lambda+1-p}
\cdot m^\Psi_K - K
\nonumber
\\
&\stackrel{\mathrm{(b)}} {\sim} & \biggl(\frac{p}{\lambda+1-p}\cdot\frac
{\lambda+1-p}{\lambda-(1-p)} -1
\biggr)\cdot K,
\nonumber
\\
&= & \frac{1-\lambda}{\lambda-(1-p)} \cdot K\qquad\mbox{a.s.}
\nonumber
\end{eqnarray}
as $K \to\infty$, where (a) follows from equations~(\ref
{eqmnbKasymp1}) and (\ref{eqtilMratio}), (b) from equations~(\ref
{eqmnbKasymp1}) and (\ref{eqsmalloK1}), which completes the proof.

\section{\texorpdfstring{Applications to resource pooling}{Applications to resource pooling}}\label{secresourcepooling}

We discuss in this section some of the implications of our results in
the context of a multi-server model for resource pooling \cite{TX12},
illustrated in Figure~\ref{figpooling}, which has partially motivated
our initial inquiry.

We briefly review the model in \cite{TX12} below, and the reader is
referred to the original paper for a more rigorous description. Fix a
coefficient $p\in[0,1]$. The system consists of $N$ stations, each of
which receives an arrival stream of jobs at rate $\lambda\in(0,1)$
and has one queue to store the unprocessed jobs. The system has a total
amount of processing capacity of $N$ jobs per unit time and is divided
between two types of servers. Each queue is equipped with a \emph
{local server} of rate $1-p$, which is capable of serving only the jobs
directed to the respective station. All stations share a \emph{central
server} of rate $pN$, which always fetches a job from the most loaded
station, following a longest-queue-first (LQF) scheduling policy. In
other words, a fraction $p$ of the total processing resources is being
\emph{pooled} in a centralized fashion, while the remainder is
distributed across individual stations. All arrival and service token
generation processes are assumed to be Poisson and independent from one
another (similarly to Section~\ref{secmodel}).

A main result of \cite{TX12} is that even a small amount of resource
pooling (small but positive $p$) can have significant benefits over a
fully distributed system ($p=0$). In particular, for any $p>0$, and in
the limit as the system size $N \to\infty$, the average delay across
the whole system scales as $\sim\log_{1/(1-p)}{\frac{1}{1-\lambda
}}$, as $\lambda\to1$; note that this is the same scaling as in
Theorem~\ref{teoonline}. This is an exponential improvement over the
scaling of ${\sim}\frac{1}{1-\lambda}$ when no resource pooling is
implemented; that is, $p=0$.

We next explain how our problem is connected to the resource pooling
model described above, and how the current paper suggests that the
results in \cite{TX12} can be extended in several directions. Consider
a similar $N$-station system as in \cite{TX12}, with the only
difference being that instead of the central server fetching jobs from
the local stations, the central server simply fetches jobs from a
``central queue,'' which stores jobs redirected from the local stations
(see Figure~\ref{figpoolingcqueue}). Denote by $ \{R_i(t)\dvtx t\in
\mathbb{R}_{+} \}$, $i\in\{1,\ldots, N\}$, the counting
process where $R_i(t)$ is the cumulative number of jobs redirected to
the central queue from station $i$ by time $t$. Assume that $\limsup_{t
\to\infty} \frac{1}{t}R_i(t) = p-\varepsilon$ almost surely for
all $i\in\{1,\ldots,N\}$, for some $\varepsilon>0$.\footnote{Since the
central server runs at rate $pN$, the rate of $R_i(t)$ cannot exceed
$p$, assuming it is the same across all $i$. }

From the perspective of the central queue, it receives an arrival
stream $R^N$, created by merging $N$ redirection streams, $R^N(t) =
\sum_{i=1}^N R_i(t)$. The process $R^N$ is of rate $(p-\varepsilon)N$,
and it is served by a service token generation process of rate $pN$.
The traffic intensity of the of central queue (arrival rate divided by
service rate) is therefore $\rho_c = (p-\varepsilon)N/pN=1-\varepsilon
/p<1$. Denote by $Q^N\in\mathbb{Z}_{+}$ the length of the central
queue in steady-state. Suppose that it can be shown that\footnote{For
an example where this is true, assume that every local station adopts a
randomized rule and redirects an incoming job to the central queue with
probability $\frac{p-\varepsilon}{\lambda}$ [and that $\lambda$ is
sufficiently close to $1$ so that$\frac{p-\varepsilon}{\lambda}\in
(0,1)$]. Then $R_i(t)$ is a Poisson process, and by the merging
property of Poisson processes, so is $R_N(t)$. This implies that the
central\vspace*{1pt} queue is essentially an $M/M/1$ queue with traffic intensity
$\rho_c = (p-\varepsilon)/p$, and we have that $\mathbb{E} (Q^N
)= \frac{\rho_c}{1-\rho_c}$ for all $N$. }
%
%
\begin{equation}
\limsup_{N\to\infty} \mathbb{E} \bigl(Q^N \bigr) <
\infty. \label{eqcentralfinite}
\end{equation}
A key consequence of equation~(\ref{eqcentralfinite}) is that, for
large values of $N$, $Q^N$ becomes negligible in the calculation of the
system's average queue length: the average queue length across the
whole system coincides with the average queue length among the \emph
{local} stations, as $N\to\infty$. In particular, this implies that,
in the limit of $N\to\infty$, the task of scheduling for the resource
pooling system could alternatively be implemented by running a separate
admissions control mechanism, with the rate of redirection equal to
$p-\varepsilon$, where all redirected jobs are sent to the central queue,
granted that the streams of redirected jobs ($R_i(t)$) are sufficiently
well behaved so that equation~(\ref{eqcentralfinite}) holds. This is
essentially the justification for the equivalence between the resource
pooling and admissions control problems, discussed at the beginning of
this paper (Section~\ref{secintroadminresourcepool}).

With this connection in mind, several implications follow readily from
the results in the current paper, two of which are given below:
\begin{longlist}[(2)]
\item[(1)] The original longest-queue-first scheduling policy employed
by the central server in \cite{TX12} is \emph{centralized}: each
fetching decision of the central server \mbox{requires} the full knowledge of
the queue lengths at all local stations. However, Theorem~\ref{teoonline} suggests that the same system-wide delay scaling in the
resource pooling scenario could also be achieved by a \emph
{distributed} implementation: each server simply runs the same
threshold policy, $\pi_{\mathrm{th}}^{L (p-\varepsilon,\lambda)}$,
and routes all deleted jobs to the central queue. To prove this
rigorously, one needs to establish the validity of equation~(\ref
{eqcentralfinite}), which we will leave as future work.

\item[(2)] A fairly tedious stochastic coupling argument was employed
in \cite{TX12} to establish a matching lower bound for the $\sim\log
_{1/(1-p)}{\frac{1}{1-\lambda}}$ delay scaling, by showing that the
performance of the LQF policy is no worse than any other online policy.
Instead of using stochastic coupling, the lower bound in Theorem~\ref{teoonline} immediately implies a lower bound for the resource pooling
problem in the limit of $N\to\infty$, if one assumes that the central
server adopts a \emph{symmetric} scheduling policy, where the it does
not distinguish between two local stations beyond their queue
lengths.\footnote{This is a natural family of policies to study, since
all local servers, with the same arrival and service rate, are indeed
identical.} To see this, note that the rate of $R_i(t)$ are identical
under any symmetric scheduling policy, which implies that it must be
less than $p$ for all $i$. Therefore, the lower bound derived for the
admissions control problem on a \emph{single queue} with a redirection
rate of $p$ automatically carries over to the resource pooling problem.
Note that, unlike the previous item, this lower bound does not rely on
the validity of equation~(\ref{eqcentralfinite}).
\end{longlist}

Both observations above exploit the equivalence of the two problems in
the regime of $N\to\infty$. With the same insight, one could also
potentially generalize the delay scaling results in \cite{TX12} to
scenarios where the arrival rates to the local stations are
nonuniform, or where future information is available. Both extensions
seem difficult to accomplish using the original framework of \cite
{TX12}, which is based on a fluid model that heavily exploits the
symmetry in the system. On the downside, however, the results in this
paper tell us very little when system size $N$ is \emph{small}, in
which case it is highly conceivable that a centralized scheduling rule,
such as the longest-queue-first policy, can out-perform a collection of
decentralized admissions control rules.
\end{appendix}

\section*{\texorpdfstring{Acknowledgment.}{Acknowledgment}}\label{sec9}
The authors are grateful for the anonymous reviewer's feedback.



\printaddresses

\end{document}